\documentclass{llncs}
\usepackage{a4wide}
\usepackage{makeidx}  
\usepackage{amsfonts}
\usepackage{amsmath}
\usepackage[dutch,USenglish]{babel}

\newcommand{\RSW}{Radcliffe, Scott and Wilmer}

\newcommand{\QeD}{$\Box$}
\newcommand{\flp}[1]{\underline{#1}}

\newif\ifcomment\commentfalse
\def\commentON{\commenttrue}

\long\outer\def\bc#1\ec{{\ifcomment \sloppy  $[${\bf suggest}]
{{#1}} \textbf{[end]} \fi }}

\long\outer\def\br#1\er{{\ifcomment \sloppy  $[${\bf suggest remove}]
{{#1}} \textbf{[end]} \fi }}

\long\outer\def\bo#1\eo{{\ifcomment \sloppy  $[${\bf instead of}]
{\textit{#1}} \textbf{[end]}  \fi }}

\long\outer\def\BC#1\EC{{\ifcomment \sloppy \par \#  \dotfill
{\textsc{#1}} \dotfill \# \par \fi }}

\commentON

\begin{document}

\title{Prefix reversals on binary and ternary strings\thanks{This research has been funded by the Dutch BSIK/BRICKS
project.}}
\titlerunning{Prefix reversals on binary and ternary strings}  
%
\author{Cor Hurkens\inst{1}, Leo van Iersel\inst{1}, Judith Keijsper\inst{1}, Steven
Kelk\inst{2}, Leen
Stougie\inst{1}\inst{2} \and John Tromp\inst{2}}

\authorrunning{Hurkens, van Iersel, Keijsper, Kelk, Stougie and Tromp}

\institute{Technische Universiteit Eindhoven (TU/e), Den Dolech 2, 5612 AX Eindhoven, Netherlands.\\
\email{wscor@win.tue.nl, l.j.j.v.iersel@tue.nl, j.c.m.keijsper@tue.nl} \\
\and
Centrum voor Wiskunde en Informatica (CWI), Kruislaan 413, 1098 SJ Amsterdam, Netherlands. \\
\email{S.M.Kelk@cwi.nl, Leen.Stougie@cwi.nl, John.Tromp@cwi.nl} \\
}

\maketitle              

\begin{abstract}
Given a permutation $\pi$, the application of prefix reversal $f^{(i)}$ to $\pi$ reverses
the order of the first $i$ elements of $\pi$.
The problem of Sorting By Prefix Reversals
(also known as \emph{pancake flipping}),
made famous by Gates and Papadimitriou
(Bounds for sorting by prefix reversal, \emph{Discrete Mathematics} 27, pp. 47-57),
asks for the minimum number of prefix reversals required to
sort the elements of a given permutation.
In this paper we study a variant of
this problem where the prefix reversals
act not on permutations but on strings
over a fixed size alphabet.
We determine the minimum number of
prefix reversals required to sort binary and ternary strings,
with polynomial-time algorithms for these sorting problems as a result;
demonstrate that
computing the minimum prefix reversal distance between two
binary strings is NP-hard; give an exact expression
for the prefix reversal diameter of binary strings,
and give bounds on the prefix reversal diameter of ternary strings.
We also consider a weaker form of sorting called \emph{grouping}
(of identical symbols) and give polynomial-time algorithms for
optimally grouping binary and ternary strings. A number of intriguing open
problems are also discussed.
\end{abstract}
\section{Introduction}
%
For a permutation $\pi=\pi(0)\pi(1)\ldots \pi(n-1)$
the application of {\em prefix reversal} $f^{(i)}$, which we call {\em flip} for short,
to $\pi$ reverses the order of the first $i$ elements:
$f^{(i)}(\pi) = \pi(i-1)\ldots \pi(0)\pi(i)\ldots \pi(n-1).$
The problem of \emph{Sorting By Prefix Reversals} (MIN-SBPR),
brought to popularity by Gates and Papadimitriou \cite{gatespap} and often referred to
as the \emph{pancake flipping problem}, is defined as follows: given a
permutation $\pi$ of $\{0,1,\ldots,n-1\}$,
determine its sorting distance i.e.
the smallest number of flips required to transform $\pi$
into the identity permutation $0 1 \ldots (n-1)$.\footnote{
We adopt the convention of numbering from 0 rather than 1.}

MIN-SBPR has practical relevance in the area of efficient network
design~\cite{heydarisud,pnetwork}, and arises
in the context of computational biology when seeking to
explain the genetic difference between two given species by the
most parsimonious (i.e. shortest) sequence of gene rearrangements.
The computational complexity of MIN-SBPR remains open.
A recent 2-approximation algorithm~\cite{2approx} is currently the
best-known approximation result\footnote{Although not
explicitly described as such, the algorithm provided ten years earlier
in~\cite{burnt} is a 2-approximation algorithm
for the \emph{signed} version of the problem.}.
Indeed, most studies to date have focused not on the computational complexity
of MIN-SBPR but rather on determining
the worst-case sorting distance $wc(n)$ over all length-$n$ permutations
i.e. the ``worst case scenario'' for length-$n$ permutations.
From \cite{gatespap} and \cite{heydarisud} we know
that $(15/14)n \leq wc(n) \leq (5n+5)/3$.

A natural variant of MIN-SBPR is to consider the action of flips
not on permutations but on strings over fixed size alphabets.
%
The shift from permutations
to strings
alters the problem universe somewhat. With permutations,
for example, the {\em distance problem}, i.e. given two permutations $\pi_1$ and
$\pi_2$, determine the smallest number of flips required to transform
$\pi_1$ into $\pi_2$, is equivalent to sorting,
because the symbols can simply be
relabelled to make either permutation equal to the identity permutation.
For strings like $101$, such a relabelling is not possible.
Thus, the distance problem on string pairs appears to be strictly
more general than the sorting problem on strings,
naturally defined as putting all elements in non-descending order.

Indeed, papers by Christie and Irving \cite{christieirving} and Radcliffe,
Scott and Wilmer \cite{rsw} explore the consequences of
switching from permutations to strings; they both
consider arbitrary (substring) reversals, and
\emph{transpositions} (where two adjacent substrings are swapped.)
It has been noted that,
viewed as a whole, such rearrangement operations on strings have bearing
on the study of orthologous gene
assignment~\cite{chen}, especially where the level of symbol repetition in the
strings is low. There is also a somewhat surprising link with the relatively
unexplored family
of \emph{string partitioning} problems~\cite{mcsp}.
To put our work in context, we briefly
describe the most relevant (for this paper) results from
\cite{christieirving} and \cite{rsw}.

The earlier paper~\cite{christieirving}, gives, in both the case of reversals
and transpositions, polynomial-time algorithms for computing
the minimum number of operations to sort a given binary string,
as well as exact, constructive diameter results on binary strings.
Additionally, their proof that computing the reversal distance between
strings is NP-hard, supports the intuition that distance problems
are harder than sorting problems on strings.
They present upper and lower bounds for computing reversal and transposition distance on binary strings.

The more recent paper~\cite{rsw} gives refined and generalised reversal
diameter results for non-fixed size alphabets. It also gives a polynomial-time
algorithm for optimally sorting a ternary (3 letter alphabet)
string with reversals. The authors refer to the prefix reversal counterparts of these (and other) results as interesting open problems. They further provide an alternative proof of Christie and Irving's NP-hardness result for reversals, and sketch a proof that computing the \emph{transposition distance} between binary strings is NP-hard. As we later note, this proof can also be used to obtain a specific reducibility result for prefix reversals. They also have some first results on
approximation (giving a PTAS - a \emph{Polynomial-Time Approximation Scheme} - for computing the distance between \emph{dense
instances}) and on the distance between random strings, both of which apply to prefix reversals as well.

In this paper we supplement results of \cite{christieirving} and
\cite{rsw} by their counterparts on prefix reversals. In Section \ref{sec:grouping} (\emph{Grouping})
we introduce a weaker form of sorting where identical symbols need only
be grouped together, while the groups can be in any order.
For grouping on binary and ternary strings we give a complete characterisation
of the minimum number of flips required to group a string, and provide polynomial-time
algorithms for computing such an optimal sequence of flips. (The complexity of grouping
over larger fixed size alphabets remains open but as an intermediate result we describe
how a PTAS can be constructed for each such problem.) Grouping aids in developing a deeper understanding of sorting
which is why we tackle it first.
It was also mentioned as a problem of interest in its own right
by Eriksson et al. \cite{bridgehand}.
Then, in Section \ref{sec:sorting}
(\emph{Sorting}), we give polynomial-time algorithms (again based on a complete
characterisation) for optimally sorting
binary and ternary strings with flips. (The complexity of sorting
also remains open for larger fixed size alphabets. As with grouping we thus provide, as an intermediate
result, a PTAS for each such problem.)
In Section \ref{sec:diam} we show that the flip diameter on binary
strings is $n-1$, and on ternary strings (for $n>3$) lies somewhere between
$n-1$ and $(4/3)n$, with empirical support for the former.
In Section \ref{sec:complexity} we show
that the flip distance problem on binary strings is NP-hard, and
point out that a reduction in \cite{rsw} also applies to prefix reversals,
showing that the flip distance problem on \emph{arbitrary} strings is
polynomial-time reducible (in an approximation-preserving sense) to
the binary problem.
We conclude in Section \ref{sec:openproblems}
with a discussion of
some of the intriguing open problems that have emerged during this work. Indeed,
our initial exploration has identified many basic (yet surprisingly difficult)
combinatorial problems that deserve further analysis.

\section{Preliminaries}
\label{sec:prelim}

Let $[k]$
denote the first $k$ non-negative integers $\{ 0,1, ..., k-1 \}$.
A $k$-ary string is a string over the alphabet $[k]$, while
a string $s$ is said to be {\em fully} $k$-ary, or to
\emph{have arity $k$}, if the set of symbols occuring in it is $[k]$.


We index the symbols in a string $s$ of length $n$ from 1 through $n$:
$s =s_1 s_2 \ldots s_n$.
Two strings are {\em compatible}
if they have the same symbol frequencies (and hence the same length),
e.g. $0012$ and $1002$ are compatible but $0012$ and $0112$ are not. For
a given string $s$, let $I(s)$ be the string obtained by sorting the
symbols of $s$ in non-descending order e.g. $I(1022011)=0011122$.
The prefix reversal (flip for short) $f^{(i)}(s)$ reverses the length $i$ prefix
of its argument, which should have length at least $i$. Alternatively,
we denote application of $f^{(i)}(s)$ by underlining the length $i$ prefix.
Thus, $f^{(2)}(2012) = \flp{20}12 = 0212$ and $f^{(3)}(2012) = \flp{201}2 = 1022$.
The \emph{flip distance} $d(s,s')$ between two strings
$s$ and $s'$ is defined as the smallest number
of flips required to transform $s$ into $s'$, if they are
compatible and $\infty$ otherwise.
Since a flip is its own inverse, flip distance is symmetric.

The \emph{flip sorting distance} $d_{\rm s}(s)=d(s ,I(s))$
of a string $s$ is defined as
the number of flips of an \emph{optimal sorting sequence} to transform $s$
into $I(s)$. An algorithm sorts $s$ optimally
if it computes an optimal sorting sequence for $s$.


In the next two sections we consider strings to be equivalent if one can be transformed into the other by repeatedly duplicating symbols and eliminating one of two adjacent identical symbols.
As representatives of the equivalence
classes we take the shortest string in each class.
These are exactly the strings in which adjacent symbols always differ.
We express all flip operations in terms of these {\em normalized}
strings. E.g. we write $f^{(3)}(2012) = \flp{201}2 = 102$.
A flip that brings two identical symbols together,
thereby shortening the string by 1, is called a {\em 1-flip},
while all others, that leave the string length invariant,
are called {\em 0-flips}.



We follow the standard notation for regular expressions:
Superindex $^i$ on a substring denotes the number of repetitions of the substring,
with  $^{*}$ and $^{+}$ denoting 0-or-more and 1-or-more repetitions, respectively,
$\epsilon$ denotes the empty string, brackets of the form $\{ \}$ are used to denote that a symbol
can be exactly one of the elements within the brackets,
and the product sign $\prod$ denotes concatenation of an indexed series.
For example $\prod_{i=1}^{3} (10^{i}2) = 102100210002$,
and $\{1,01\}^{*}\{\epsilon,0\}$ denotes the
set of binary strings with no 00 substring.

\section{Grouping}
\label{sec:grouping}

The task of sorting a string can be broken down
into two subproblems:
{\em grouping} identical symbols together and
putting the groups of identical symbols in the right {\em order}.
Notice that first grouping and then ordering may not be the most efficient way to sort strings.
Although grouping appears to be slightly easier than the sorting problem,
essentially the same questions remain open as in sorting.
Grouping binary strings is trivial
and in Section \ref{subsec:ternarygrouping} we
give the grouping distances of all ternary strings.
As a result we give polynomial time algorithms for binary and ternary grouping.
For larger alphabets the grouping problem remains open; as an intermediate result
we describe in Section
\ref{subsec:largeralphabets} a PTAS for each such problem. While the problems of grouping
and sorting are closely related for strings on
small alphabets, the problems diverge when alphabet size approaches
the string length, with permutations being the limit.

Recall that we consider only normalized strings, as representatives
of equivalence classes.
The \emph{flip grouping distance} $d_{\rm g}(s)$
of a fully $k$-ary string $s$ is defined as the minimum number
of flips required to reduce the string to one of length $k$.

\subsection{Grouping binary and ternary strings}
\label{subsec:ternarygrouping}

\begin{lemma}
\label{grouping:lem:lowerbound} $d_{\rm g}(s)\geq n-k$
for any fully $k$-ary string $s$ of length $n$.
\end{lemma}
\begin{proof}
The proof follows from the observations that, after grouping, fully $k$-ary string $s$
has length $k$ and that each flip can shorten $s$ by at most 1.
\QeD
\end{proof}
\begin{lemma}
\label{grouping:lem:greedy} $d_{\rm g}(s)\leq n-2$
for any fully $k$-ary string $s$ of length $n$.
\end{lemma}
\begin{proof}
Consider the following simple algorithm. If the leading symbol occurs elsewhere then
a 1-flip bringing them together exists, so perform this 1-flip. If not, then we use a 0-flip
to put this symbol in front of a suffix in which we accumulate uniquely
appearing symbols. Repeat until the string is grouped.

Clearly no more than $n-k$ 1-flips will be necessary. Also, no
more than $k-2$ 0-flips will ever be necessary, because after $k-2$
0-flips the prefix of the string will consist of only two types of symbol, and the algorithm will never
perform a 0-move on such a string. Thus at most $(n-k) + (k-2) = n-2$
flips in total will be needed.
\QeD
\end{proof}
As a corollary we obtain the grouping distance of binary strings.
\begin{theorem}
$d_{\rm g}(s)=n-2$
for any fully binary string $s$ of length $n$. \QeD
\end{theorem}
We will now define a class of bad ternary strings and
prove that these are the only ternary strings that need $n-2$
rather than $n-3$ flips to be grouped.
\begin{definition}
\label{grouping:def:badstrings} We define \emph{bad} strings as all fully ternary strings of one of the following types, up to relabeling:
\begin{enumerate}
\item [I.] strings of length greater than 3,
  in which the leading symbol appears only once:
$0(12)^{\geq 2}$ and $02(12)^+$
\item [II.] strings having identical symbols at every other position,
starting from the last:
$(\{0,1\}2)^+$ and $(2 \{0,1\})^+ 2$
\item [III.] odd length strings whose leading symbol appears exactly once more,
at an even position, and both occurrences are followed by the same symbol:
$0(21)^+02(12)^*$
\item [IV.] the following strings:\\
$X_1 = 210212$, $X_2 = 021012$, $X_3 = 0120212$, $X_4 = 1201212$,
$X_5 = 02101212$, $X_6 = 20210212$, $X_7 = 020210212$, $X_8 = 120120212$.
\end{enumerate}
All other fully ternary strings are {\em good}.
Strings of type I, II and III, shortly I-, II-, and III-strings, respectively,
are called generically bad, or {\em g-bad} for short.
\end{definition}

\begin{lemma}
\label{grouping:lem:badstrings}
$d_{\rm g}(s)=n-2$
if ternary string $s$ of length $n$ is bad.
\end{lemma}
\begin{proof}
Because of Lemmas~\ref{grouping:lem:lowerbound} and
\ref{grouping:lem:greedy}, it suffices to show that
in each case a 0-flip is necessary:
I-strings admit only 0-flips.
A 1-flip on a II-string leads
to a II-string and eventually to a I-string.
Any III-string admits only one 1-flip
leading to a II-string.
For IV-strings, Table~\ref{grouping:tab:goodmoves}
shows that each possible 1-flip leads to either
a shorter IV-string, or to a I-,II-, or III-string.
\QeD
\end{proof}

\begin{table}[h]
\begin{centering}
\begin{tabular}{|l|l|}
\hline
$X_1$ & $X_6$\\ \hline
$\flp{210}212=01212$ is of type I & $\flp{20}210212 = 0210212$ is of type III\\
$\flp{21021}2 = 12012$ is of type III & $\flp{20210}212 = 0120212 = X_3$\\
\cline{1-1}
$X_2$ & $\flp{2021021}2 = 1201202$ is of type III\\ \hline
$\flp{021}012 = 12012$ is of type III & $X_7$\\ \hline
$X_3$ & $\flp{02}0210212 = 20210212 = X_6$\\ \cline{1-1}
$\flp{012}0212 = 210212 = X_1$ & $\flp{02021}0212 = 12020212$ is of type II\\ \hline
$X_4$ & $X_8$\\ \hline
$\flp{120}1212 = 021212$ is of type I and II\ \ \ \ \  & $\flp{120}120212 = 02120212$ is of type II\ \ \ \ \ \\
$\flp{12012}12 = 210212 = X_1$ & $\flp{1201202}12 = 20210212 = X_6$\\ \cline{1-1}
$X_5$ & \\ \cline{1-1}
$\flp{021}01212 = 1201212 = X_4$ & \\
\hline
\end{tabular}\\
\end{centering}
\medskip

\caption{type IV strings and all their 1-flips.}
\label{grouping:tab:goodmoves}
\end{table}

\newpage

\begin{lemma}
\label{grouping:lem:goodstrings}
$d_{\rm g}(s)=n-3$
if ternary string $s$ of length $n$ is good.
\end{lemma}
\begin{proof}
The proof is by induction on $n$. The induction basis for $n=3$ is
trivial.
We show the statement for strings of length $n+1$ by showing that if
a bad string $s'$ of length $n$ can be obtained through a 1-flip from
a good ({\em parent}) string $s$ of length $n+1$,
then $s$ admits another 1-flip which leads to
a good string.
Note that a 1-flip $f^{(i)}(s)=s'$ brings symbols $s_1$ and $s_{i+1}$ together,
hence $s_1 = s_{i+1} \neq s_i = s'_1$ which shows that the symbol
deleted from {\em parent} $s$ differs from the leading symbol of {\em child} $s'$.
We enumerate all possible bad child strings $s'$ and distinguish cases based on
the leading symbol of good parent $s$.

For IV-strings,
Table~\ref{grouping:excep} lists all parents with, for each good parent,
a 1-flip to a good string.
It remains to prove that for each g-bad string all parents are either bad
or have a g-1-flip, defined as
a 1-flip resulting in a string that is not g-bad
(i.e. either good or of type IV).
\begin{description}

\item [Type I, odd:] $0(12)^{\geq 2}$ has possible parents starting with:
  \begin{enumerate}
  \item[1:] $1(21)^i012(12)^j$ with $i+j>0$:
     \begin{description}
     \item[If] $i>0$ there is a g-1-flip
    $\flp{12}1(21)^{i-1}012(12)^j = (21)^i012(12)^j$;
     \item[If] $i=0$ and $j>0$ there is a g-1-flip
    $\flp{1012}(12)^j = 210(12)^j$;
     \end{description}
  \item[2:] $21(21)^i02(12)^j$ with $i+j>0$.
     \begin{description}
     \item[If] $i>0$  there is a g-1-flip
    $\flp{21}(21)^i02(12)^j = 1(21)^i02(12)^j$;
     \item[If] $i=0$ and $j>1$ there is a g-1-flip
    $\flp{21021}2(12)^{j-1} = 120(12)^j$;
     \item[If] $i=0$ and $j=1$ the parent is $210212 = X_1$.
     \end{description}
  \end{enumerate}

\item [Type I, even:] these strings are also of type II, see below.

\item [Type II, odd:] $(2\{0,1\})^+2$
    has only parents of type II.

\item [Type II, even:] $02(\{0,1\}2)^*$ has possible parents starting with:
  \begin{enumerate}
  \item [2:] $2(\{0,1\}2)^*$ is of type II;
  \item [1:] $12(\{0,1\}2)^*012(\{0,1\}2)^*$ with three cases for a possible third 1:
     \begin{description}
     \item[None:] parent is $12(02)^*012(02)^*$, which is of type III;
     \item[Before $01$:] then there is a g-1-flip \\
    $\flp{12(\{0,1\}2)^*}12(\{0,1\}2)^*012(\{0,1\}2)^* = 2(\{0,1\}2)^*12(\{0,1\}2)^*012(\{0,1\}2)^*$;
     \item[After $01$:] then there is a g-1-flip \\
    $\flp{12(\{0,1\}2)^*012(\{0,1\}2)^*}12(\{0,1\}2)^* = 2(\{0,1\}2)^*102(\{0,1\}2)^*12(\{0,1\}2)^*$.
     \end{description}
  \end{enumerate}

\item [Type III:] $0(21)^+02(12)^*$ has possible parents starting with:
  \begin{enumerate}
  \item[1:] $(12)^i01(21)^j02(12)^k$ with $i>0$:
      \begin{description}
      \item [If] $i>1$ there is a g-1-flip
   $\flp{12}(12)^{i-1}01(21)^j02(12)^k=2(12)^{i-1}01 (21)^j02(12)^k$;
      \item [If] $i=1,j>0$ there is a g-1-flip
   $\flp{12012}1(21)^{j-1}02(12)^k=21021(21)^{j-1}02(12)^k$;
      \item [If] $i=1,j=0,k>0$ there is a g-1-flip
   $\flp{120102}(12)^k=20102(12)^k$;
      \item [If] $i=1,j=k=0$ then the parent is $120102=X_2$ (relabelled);
      \end{description}

   \item[1:] $(12)^+0(12)^+0(12)^+$: there is a g-1-flip
     $\flp{(12)^+0}(12)^+0(12)^+ = 0(21)^+20(12)^+$;

   \item[2:] $2(12)^*0(21)^+02(12)^*$: there is a g-1-flip
       $\flp{2(12)^*0(21)^+0}2(12)^* = 0(12)^+0(21)^*2$;

   \item[2:] $(21)^i20(12)^j02(12)^k$ with $j>0$:
     \begin{description}
     \item [If] $i=0, j=1$ then the parent is $210212 = X_1$;
     \item [If] $i+j>1$ then
      $\flp{(21)^i201}2(12)^{j-1}02(12)^k = 102(12)^{i+j-1}02(12)^k$ is a g-1-flip. \QeD
     \end{description}
   \end{enumerate}
\end{description}
\end{proof}

\begin{table}[h]
\begin{centering}
\begin{tabular}{|l|l||l|l||l|l|}
\hline
$X_1$    & Parents       & $X_4$     & Parents        & $X_7$       & Parents\\
\hline
$210212$ & $\flp{12102}12$    & $1201212$ & $\flp{21201}212$    & $020210212$ & $\flp{2020210}212$\\
         & $0120212=X_3$ &           & $02101212=X_5$ &             & $\flp{2020210}212$\\
         & $1201212=X_4$ &           & $\flp{2102121}2$    &             & $\flp{12020102}12$\\
\cline{1-2}
$X_2$    & Parents       &           & $\flp{21}210212$    &             & $\flp{201}2020212$\\
\cline{1-4}
$021012$ & $\flp{202101}2$    & $X_5$     & Parents        &             & $\flp{120}1202012$\\
\cline{3-4}
         & $\flp{12010}12$    & $02101212$& $\flp{202101}212$   &             & $\flp{20210}21212$\\
\cline{5-6}
         & $\flp{10}12012$    &           & $\flp{12010}1212$   & $X_8$       & Parents      \\
\cline{5-6}
         & $\flp{2010}212$    &           & $\flp{10}1201212$   & $120120212$ & $\flp{21201}20212$\\
\cline{1-2}
$X_3$    & Parents       &           & $\flp{2101}20212$   &             & $\flp{021012}0212$\\
\cline{1-2}
$0120212$& $\flp{101202}12$   &           & $\flp{12}1012012$   &             & $\flp{210}2120212$\\
         & $\flp{21020}212$   &           & $\flp{20}2010212$   &             & $\flp{021}0210212$\\
\cline{3-4}
         & $20210212=X_6$& $X_6$     & Parents        &             & $\flp{20210}21212$\\
\cline{3-4}
         & $\flp{1202}1012$   & $20210212$& $020210212=X_7$&             & $\flp{21}20210212$\\
         & $\flp{20}212012$   &           & $\flp{1202102}12$   &             &              \\
         &               &           & $\flp{012}020212$   &             &              \\
         &               &           & $120120212=X_8$&             &              \\
 \hline
\end{tabular}\\
\end{centering}
\medskip

\caption{Type IV strings, their parents, and for each good parent,
a 1-flip to a good string.} \label{grouping:excep}
\end{table}


\noindent
The following theorem results directly from the above lemmas.
\begin{theorem}
\label{grouping:thm:alg} $d_{\rm g}(s)=n-2$
if and only if fully ternary string $s$ of length $n$
is bad and $d_{\rm g}(s)=n-3$ otherwise. Moreover, there exists a
polynomial time algorithm for
grouping ternary strings with a minimum number of flips.
\end{theorem}
\begin{proof}
The first statement is direct from Lemmas~\ref{grouping:lem:badstrings}
and \ref{grouping:lem:goodstrings}. In case string $s$ is bad,
which by Definition \ref{grouping:def:badstrings} can be decided in
polynomial time, the algorithm implicit in the proof of
Lemma~\ref{grouping:lem:greedy} shows how to group $s$ optimally in polynomial time.
Otherwise, we repeatedly find a 1-flip to a good string
as guaranteed by Lemma~\ref{grouping:lem:goodstrings}.
The time complexity is $O(n^3)$, since grouping distance,
number of choices for a 1-flip, and time to perform a flip and
test whether its result is good are all $O(n)$.
\QeD
\end{proof}

\subsection{Grouping strings over larger alphabets}
\label{subsec:largeralphabets}
Lemmas~\ref{grouping:lem:lowerbound} and \ref{grouping:lem:greedy}
say that $n-k \leq d_{\rm g}(s) \leq n-2$ for any fully $k$-ary string $s$.
For any $k$ there are fully $k$-ary strings that have flip grouping
distance equal to $n-2$. For example the length $n=2(k-1)$ string
$1020\ldots(k-1)0$ requires for every 1-flip to bring a 0 to
the front first and hence we need as many 0-flips
as 1-flips, and $d_{\rm g}(1020\ldots(k-1)0)\geq 2(k-2)=2k-4=n-2$. Computer calculations
suggest that for $k=4$ and $k=5$, for $n$ large enough,
the strings with grouping distance $n-2$ are precisely those having identical symbols at every other position, starting from the last
(i.e. type II of Definition \ref{grouping:def:badstrings}).
Proving (or disproving) this statement remains open,
as well as finding a polynomial time algorithm for grouping $k$-ary strings for
any fixed $k>3$. We do, however, have the following intermediate result:

\begin{theorem}
\label{grouping:thm:ptas}
For every fixed $k$ there is a PTAS for grouping $k$-ary strings.
\end{theorem}
\begin{proof}
We show that, for every fixed $k$ and for every fixed $\epsilon > 0$ there is a
polynomial-time algorithm that, given any $k$-ary string $s$ of length $n$, computes a sequence of flips
which groups $s$ in at most $(1+\epsilon)d_{\rm g}(s)$ flips.
We assume $k\geq 4$ because for $k=2$ and $k=3$ the exact algorithms suffice.
Let $N = (k-2)/\epsilon + k$. We distinguish two cases.\\
\\
\textbf{Case 1.}
If $n \geq N$ we use the simple, ``greedy'' algorithm described in the proof of Lemma
\ref{grouping:lem:greedy}. This will group $s$ in $d_{\rm g}^G(s)$ flips with
$d_{\rm g}^G \leq n-2$ steps. This together with the lower bound of
$n-k$ on $d_{\rm g}(s)$ from Lemma \ref{grouping:lem:lowerbound}
gives $d_{\rm g}^G(s)\leq d_{\rm g}(s) + (k-2) \leq
(1+\epsilon)d_{\rm g}(s).$\\
\\
\textbf{Case 2.} If $n<N$ we compute $d_{\rm g}(s)$ by
a brute force algorithm which
simply chooses the best amongst all possible flip sequences of length $n-2$: there are $n^{n-2}$ of these.
This yields the optimal solution since $d_{\rm g}(s) \leq n-2$ (Lemma
\ref{grouping:lem:greedy}). The running time in this case is bounded by a constant.
\QeD
\end{proof}

Clearly, there is a strong relationship between grouping and sorting.
Understanding grouping may help us to understand sorting, and
lead to improved bounds (especially as the length of strings becomes large
relative to their arity), because for a $k$-ary string $s$, we have
$d_{\rm g}(s) \leq d_{\rm s}(s) \leq d_{\rm g}(s) + wc(k)$, with
$wc(k)$ the flip diameter on permutations
with $k$ elements, as defined before.

Also $d_{\rm g}(s) = \min \{ d_{\rm s}(t): t \mbox{ a relabeling of } s \}$,
which gives (for fixed $k$) a polynomial time reduction from grouping to
sorting. Thus every polynomial time algorithm for sorting by prefix
reversals directly gives a polynomial time algorithm for the grouping problem
(for fixed $k$).

\section{Sorting}
\label{sec:sorting}

In this section we present results on sorting similar to those on grouping
in the previous section. Also flip sorting distance
remains open for strings over alphabets of size larger than 3. As an
intermediate result we thus provide at the end of this section a PTAS for each
such problem.

Again a {\em 1-flip} brings identical symbols together and thus shortens the representative of the equivalence class
under symbol duplication. But since symbol order matters for sorting, relabelled strings are no longer equivalent. As
in grouping, sorting of binary strings is straightforward:
\begin{theorem}
\label{sorting:thm:bin}  $d_{\rm s}(s)=n-2$ for every fully binary string $s$ of length $n$ with $s_n=1$, and $d_{\rm s}(s)=n-1$ otherwise.
\end{theorem}
\begin{proof}
Exactly $n-2$ 1-flips suffice and are
necessary to arrive at length 2 string 01 or 10.
If the last symbol is 0 an additional 0-flip is necessary
putting a 1 at the end. All these flips can be $f^{(2)}$.
\QeD
\end{proof}

From Lemma~\ref{grouping:lem:lowerbound} we know that
$d_{\rm g}(s)\geq n-3$ and hence $d_{\rm s}(s)\geq n-3$
for every ternary string $s$ of length $n$.
In the upper bound on $d_{\rm s}(s)$ we derive below
we focus on strings $s$ ending in
a $2$ ($s_n=2$), since sorting distance is invariant under appending a $2$ to a string.
It turns out that, when sorting a ternary
string ending in a $2$, one needs at most one 0-flip, except for the string $0212$.
\begin{lemma}
\label{sorting:lem:upperbound} $d_{\rm s}(s)\leq n-2$ for every fully ternary string $s$
of length $n$ with $s_n=2$, except $0212$.
\end{lemma}
\begin{proof} It is easy to check that $0212$ requires 3 flips to be sorted.
By induction on $n$ we prove the rest of the lemma. The basis case of $n=3$ is trivial.
For a string $s$ of length $n>3$ we distinguish three cases:
\begin{itemize}
\item $s_{n-1}=0$:
If $s=20102$ it is sorted in 3 flips:
$\flp{2010}2 \rightarrow \flp{01}02 \rightarrow \flp{10}2 \rightarrow 012$.
Otherwise, by induction and relabeling $0\leftrightarrow 2$,
the string $s_1\ldots s_{n-1}$ can be reduced to $210$ in $n-3$ flips
(to 20 or 10 by Theorem~\ref{sorting:thm:bin} if $s_1\ldots s_{n-1}$
has only two symbols), and one more flip sorts $s$ to $012$.
\item $s_{n-1}=1$, $s_1=0$ and appears only once:
Thus $s=0(12)^{\geq 2}$ or $s=02(12)^{\geq 2}$. Then
$s$ can be sorted with only one 0-flip:
$\flp{0(12)^+1}2 \rightarrow \flp{1(2}1)^+02 \rightarrow \ldots \rightarrow
\flp{210}2 \rightarrow 012 $ or, respectively,
    $\flp{02}(12)^{\geq 2} \rightarrow \flp{20(12)^+1}2
    \rightarrow \flp{(12)^+}102 \rightarrow \ldots \rightarrow \flp{210}2
\rightarrow 012$.

\item $s_{n-1}=1$, $s_1$ not unique: \\
    If $s=12012$ then 3 flips suffice:
    $\flp{12}012 \rightarrow \flp{2101}2 \rightarrow
     \flp{10}12 \rightarrow 012$. \\
    Otherwise, since the other 2 parents of $0212$ can flip to $1202$,
    there is a 1-flip to a string $\neq 0212$ to which
    we can apply the  induction hypothesis. \QeD
\end{itemize}
\end{proof}

As in Section \ref{sec:grouping}, we characterise the strings ending in a $2$
that need $n-2$ rather than $n-3$ flips to sort.

\begin{definition}
\label{sorting:def:badstrings}
We define {\em bad} strings as all fully ternary strings ending in a $2$
of the types:
\begin{enumerate}
\item [I.] $0(12)^{\geq 2}$
\item [II.] $(\{0,1\}2)^+$ and $2(\{0,1\}2)^+$
\item [III.] $(\{1,2\}0)^+2$ and  $0(\{1,2\}0)^+$
\item [IV.] $(\{1,2\}0)^+12$ and $(0\{1,2\})^+012$ with at least two $2$s.
\item [V.] $(01)^*0212$  and   $(10)^+212$
\item [VI.] $1(20)^+1(20)^*2$  and $0(21)^+0(21)^*2$
\item [VII.] $1(02)^+1(02)^+$
\item [VIII.] $1(02)^+12$
\item [IX.] 77 strings of length at most 11, shown in Table \ref{sorting:excep}.
\end{enumerate}
All other fully ternary strings ending in a $2$ are {\em good} strings. Strings of type I-VIII (I-strings ...
VIII-strings for short) are called generically bad, or {\em g-bad} for short.
\end{definition}

\begin{table}[h]
\begin{centering}
\begin{tabular}{|l|l|l|l|}
\hline
$   Y_{ 1   }=  210212  $&$ Y_{ 21  }=   10212012 $&$ Y_{ 41  }=  021202012    $&$ Y_{ 61  }=  0210212012   $\\
$   Y_{ 2   }=  021012   $&$ Y_{ 22  }=  02121012 $&$ Y_{ 42  }=  021201012    $&$ Y_{ 62  }=  1021202012  $\\
$   Y_{ 3   }=  212012  $&$ Y_{ 23  }=   02120102 $&$ Y_{ 43  }=  020210212    $&$ Y_{ 63  }=  1021201012  $\\
$   Y_{ 4   }=  120102  $&$ Y_{ 24  }=   10102102 $&$ Y_{ 44  }=  101020212    $&$ Y_{ 64  }=  1020210212  $\\
$   Y_{ 5   }=  201202  $&$ Y_{ 25  }=   02010212 $&$ Y_{ 45  }=  020212012    $&$ Y_{ 65  }=  1010210202  $\\
$   Y_{ 6   }=  0210202  $&$ Y_{ 26  }=  21202012 $&$ Y_{ 46  }=  212010202    $&$ Y_{ 66  }=  0202010212   $\\
$   Y_{ 7   }=  1021202 $&$ Y_{ 27  }=   21201012 $&$ Y_{ 47  }=  212012012    $&$ Y_{ 67  }=  2120202012  $\\
$   Y_{ 8   }=  0212012  $&$ Y_{ 28  }=  21201202 $&$ Y_{ 48  }=  010210212    $&$ Y_{ 68  }=  2120102012  $\\
$   Y_{ 9   }=  2120102 $&$ Y_{ 29  }=   20210212 $&$ Y_{ 49  }=  010210202    $&$ Y_{ 69  }=  2021021212  $\\
$   Y_{ 10  }=  0102102  $&$ Y_{ 30  }=  01021202 $&$ Y_{ 50  }=  010212012    $&$ Y_{ 70  }=  2010212012  $\\
$   Y_{ 11  }=  1212012 $&$ Y_{ 31  }=   01020212 $&$ Y_{ 51  }=  202010212    $&$ Y_{ 71  }=  1201021202  $\\
$   Y_{ 12  }=  2010212 $&$ Y_{ 32  }=   20212012 $&$ Y_{ 52  }=  121202012    $&$ Y_{ 72  }=  1201202012  $\\
$   Y_{ 13  }=  0120212  $&$ Y_{ 33  }=  12120102 $&$ Y_{ 53  }=  121201202    $&$ Y_{ 73  }=  10202010212 $\\
$   Y_{ 14  }=  1201012 $&$ Y_{ 34  }=   12010212 $&$ Y_{ 54  }=  201021202    $&$ Y_{ 74  }=  02120102012  $\\
$   Y_{ 15  }=  1201212 $&$ Y_{ 35  }=   12010202 $&$ Y_{ 55  }=  120212012    $&$ Y_{ 75  }=  02021021212  $\\
$   Y_{ 16  }=  2012012 $&$ Y_{ 36  }=   20120102 $&$ Y_{ 56  }=  012021212    $&$ Y_{ 76  }=  21201202012 $\\
$   Y_{ 17  }=  10210212 $&$ Y_{ 37  }=  12012012 $&$ Y_{ 57  }=  120102012    $&$ Y_{ 77  }=  12120202012 $\\
$   Y_{ 18  }=  21021212 $&$ Y_{ 38  }=  021021202$&$ Y_{ 58  }=  201202012    $&$                 $\\
$   Y_{ 19  }=  02102012 $&$ Y_{ 39  }=  102120102$&$ Y_{ 59  }=  120120212    $&$                 $\\
$   Y_{ 20  }=  02101212 $&$ Y_{ 40  }=  102010212$&$ Y_{ 60  }=  201201012    $&$                 $\\
\hline
\end{tabular}\\
\end{centering}
\medskip
\caption{Type IX strings}
\label{sorting:excep}
\end{table}

This definition makes $0212$ a bad string as well.
From Lemma \ref{sorting:lem:upperbound} we know
that $0212$ is the only ternary string
ending in a $2$ with sorting distance $n-1$.

\begin{theorem}
\label{thm:diamref}
String $0212$ has sorting distance $3$.
Any other fully ternary string $s$ of length $n$ with $s_n=2$
has prefix reversal sorting distance $n-2$ if it is bad and
$n-3$ if it is good.
A fully ternary string $s$ ending in a $0$ or $1$
has the same sorting distance as $s2$.
\end{theorem}
\begin{proof}
Directly from Lemmas~\ref{sorting:lem:sortbad}
and \ref{sorting:lem:sortgood} below.
Note that every sorting sequence for
$s$ sorts $s2$ as well while every sorting sequence for $s2$
can be modified to avoid flipping the whole string and thus works for $s$ as well.
\QeD
\end{proof}

\begin{lemma}
\label{sorting:lem:sortbad} $d_{\rm s}(s)=n-2$
for every bad ternary string $s\neq 0212$ of length $n$.
\end{lemma}
\begin{proof}
Since $d_{\rm s}(s)\geq n-3$ and any 1-flip decreases the length of the string by $1$, Lemma~\ref{sorting:lem:upperbound}
says it suffices to show that for each type in Definition \ref{sorting:def:badstrings}
a 0-flip is necessary.
\begin{itemize}
\item For I-strings only 0-flips are possible.

\item A 1-flip on a II- or III-string leads to a string of the same type, so that eventually no 1-flip is possible.

\item A 1-flip on a IV-string leads either again to a IV-string or (when destroying the $12$ suffix) to a III-string.

\item A 1-flip on a V-string leads either again to a V-string or (when destroying the suffix with a $\flp{\ldots
02}12$ flip) to a IV-string. Flips $\flp{\ldots 0}212$ and $\flp{\ldots 021}2$ are not possible for lack of more
$2$'s.

\item For strings of VI-, VII- and VIII-strings only one 1-flip is possible, leading to II-, III- and IV-strings
respectively.

\item For IX-strings, Table \ref{sorting:veriftab} in the appendix lists all possible 1-flips ultimately leading to a
string of type I-VIII. \QeD
\end{itemize}
\end{proof}

\begin{lemma}
\label{sorting:lem:sortgood} $d_{\rm s}(s)=n-3$
for every good ternary string $s$ of length $n$.
\end{lemma}
\begin{proof}
The proof is by induction on $n$ and is similar to the proof of Lemma \ref{grouping:lem:goodstrings}. The induction
basis for $n=3$ is again trivial. We prove that for each g-bad string of length $n$ all parents (of length $n+1$) are
either bad or have a 1-flip to a string that is not g-bad (i.e. either good or of type IX). Remember that such a flip
is called a \emph{g-1-flip}. That for each IX-string all parents are either bad or have a 1-flip to a good string is
proved by case checking in Table~\ref{sorting:veriftab} (see appendix). Together this proves that every good string of
length $n+1$ has a 1-flip to a good string of length $n$ and therefore the lemma.
\begin{description}
\item [Type I:] $0(12)^+$ has possible parents starting with:
    \begin{enumerate}
    \item [1:] $1(21)^i0(12)^j$ with $j>0$:
        \begin{description}
        \item [If] $i>0$: there is a g-1-flip $\flp{12}1(21)^{i-1}0(12)^j=(21)^i0(12)^j$;
        \item [If] $i=0$, $j>1$: there is a g-1-flip $\flp{1012}(12)^{j-1}=210(12)^{j-1}$;
        \item [If] $i=0$, $j=1$: there is a g-1-flip $\flp{10}12=012$;
        \end{description}
    \item [2:] $(21)^i02(12)^j$ with $i>0$:
        \begin{description}
        \item [If] $i>1$: there is a g-1-flip $\flp{21}21(21)^{i-2}02(12)^j=1(21)^{i-1}02(12)^j$;
        \item [If] $i=1$, $j>0$: there is a g-1-flip $\flp{21021}2(12)^{j-1}=120(12)^j$;
        \item [If] $i=1$, $j=0$: there is a g-1-flip $\flp{210}2=012$.
        \end{description}
    \end{enumerate}
\item [Type II, even:] $(\{0,1\}2)^+$ has possible parents starting with:
    \begin{enumerate}
    \item [0:] $0(2\{0,1\})^*2102(\{0,1\}2)^*$, with three cases for a possible third $0$:
        \begin{description}
        \item [None:] the parent is of type VI;
        \item [Before $2102$:] there is a g-1-flip\\
        $\flp{0(2\{0,1\})^*2}0(2\{0,1\})^*2102(\{0,1\}2)^*=2(\{0,1\}2)^*0(2\{0,1\})^*2102(\{0,1\}2)^*$;
        \item [After $2102$:] there is a g-1-flip\\
        $\flp{0(2\{0,1\})^*2102(\{0,1\}2)^*}02(\{0,1\}2)^*=(2\{0,1\})^*2012(\{0,1\}2)^*02(\{0,1\}2)^*$;
        \end{description}
    \item [1:] $1(2\{0,1\})^*2012(\{0,1\}2)^*$, with three cases for a possible third $1$:
        \begin{description}
        \item [None:] the parent is of type VI;
        \item [Before $2012$:] there is a g-1-flip\\
        $\flp{1(2\{0,1\})^*2}1(2\{0,1\})^*2012(\{0,1\}2)^*=2(\{0,1\}2)^*1(2\{0,1\})^*2012(\{0,1\}2)^*$;
        \item [After $2012$:] there is a g-1-flip\\
        $\flp{1(2\{0,1\})^*2012(\{0,1\}2)^*}12(\{0,1\}2)^*=(2\{0,1\})^*2102(\{0,1\}2)^*12(\{0,1\}2)^*$;
        \end{description}
    \item [2:] $2\{0,1\}(2\{0,1\})^*2(\{0,1\}2)^*$ is of type II.
    \end{enumerate}
\item [Type II, odd:] $2(\{0,1\}2)^+$ has possible parents starting with:
    \begin{enumerate}
    \item [0:] $0(2\{0,1\})^*202(\{0,1\}2)^*$ is of type II;
    \item [1:] $1(2\{0,1\})^*212(\{0,1\}2)^*$ is of type II.
    \end{enumerate}
\item [Type III, even:] $0(\{1,2\}0)^+2$ has possible parents starting with:
    \begin{enumerate}
    \item [1:] $1(0\{1,2\})^+010(\{1,2\}0)^*2$ is of type III;
    \item [2:] $2(0\{1,2\})^+020(\{1,2\}0)^*2$ is of type III;
    \item [2:] $2(0\{1,2\})^+02$ is of type III.
    \end{enumerate}

\item [Type III, odd:] $(\{1,2\}0)^+2$ has possible parents starting with:
    \begin{enumerate}
    \item [0:] $0\{1,2\}(0\{1,2\})^*0(\{1,2\}0)^*2$ is of type III;
    \item [1:] $1(0\{1,2\})^*0210(\{1,2\}0)^*2$, there are three cases for a possible third $1$:
        \begin{description}
        \item [None:] the parent is of type VII;
        \item [Before $0210$:] there is a g-1-flip\\
        $\flp{1(0\{1,2\})^*0}1(0\{1,2\})^*0210(\{1,2\}0)^*2=0(\{1,2\}0)^*1(0\{1,2\})^*0210(\{1,2\}0)^*2$;
        \item [After $0210$:] there is a g-1-flip\\
        $\flp{1(0\{1,2\})^*0210(\{1,2\}0)^*}10(\{1,2\}0)^*2=(0\{1,2\})^*0120(\{1,2\}0)^*10(\{1,2\}0)^*2$;
        \end{description}
    \item [2:] $\flp{2(0\{1,2\})^*0120(\{1,2\}0)^*}2=(0\{1,2\})^*0210(\{1,2\}0)^*2$ is a g-1-flip\\
    unless this last string is $02102$ (type VI), but then the parent is $201202=Y_5$;
    \item [2:] $2(0\{1,2\})^*012$ is of type IV.
    \end{enumerate}

\item [Type IV, even:] $(\{1,2\}0)^+12$ with a second $2$, has possible parents starting with:
    \begin{enumerate}
    \item [0:] $0({1,2}0)^+12$, with a second $2$, is of type IV;
    \item [1:] $\flp{1(0\{1,2\})^*0210(\{1,2\}0)^*}12=(0\{1,2\})^*0120(\{1,2\}0)^*12$ is a g-1-flip;
    \item [1:] $1(0\{1,2\})^*0212$, with three cases:
        \begin{description}
        \item [No third $2$:] the parent is of type V;
        \item [No third $1$:] the parent is of type VIII;
        \item [Otherwise:] $\flp{1(0\{1,2\})^*0}1(0\{1,2\})^*0212=0(\{1,2\}0)^*1(0\{1,2\})^*0212$ (with a third $2$) is a g-1-flip;
        \end{description}
    \item [2:] $2(0\{1,2\})^*0120(\{1,2\}0)^*12$, with four cases:
        \begin{description}
        \item [A fourth 2 before $0120$:] there is a g-1-flip\\
        $\flp{2(0\{1,2\})^*0}2(0\{1,2\})^*0120(\{1,2\}0)^*12=0(\{1,2\}0)^+120(\{1,2\}0)^*12$;
        \item [A fourth 2 after $0120$:] there is a g-1-flip\\
        $\flp{2(0\{1,2\})^*0120(\{1,2\}0)^*}20(\{1,2\}0)^*12=(0\{1,2\})^*0210(\{1,2\}0)^+12$;
        \item [A third 1:] $\flp{2(0\{1,2\})^*0120(\{1,2\}0)^*1}2=1(0\{1,2\})^*0210(\{1,2\}0)^*2$ is a g-1-flip;
        \item [Otherwise:] $2012012=Y_{16}$;
        \end{description}
    \item [2:] $\flp{21(0\{1,2\})^*0}2(0\{1,2\})^*012=0(\{1,2\}0)^*12(0\{1,2\})^*012$ is a g-1-flip.
    \end{enumerate}

\item [Type IV, odd:] $0(\{1,2\}0)^+12$ with a second $2$, has possible parents starting with:
    \begin{enumerate}
    \item [1:] $1(0\{1,2\})^+012$, with a second 2, is of type IV;
    \item [2:] $2(0\{1,2\})^+012$ is of type IV;
    \item [2:] $\flp{21(0\{1,2\})^*0}2(0\{1,2\})^*02=0(\{1,2\}0)^*12(0\{1,2\})^*02$ is a g-1-flip.
    \end{enumerate}

\item [Type V, even:] $0(10)^+212$ ($0212$ is also of type I), has possible parents starting with:
    \begin{enumerate}
    \item [1:] $(10)^+212$ is of type V;
    \item [1:] $\flp{120}1(01)^*012=021(01)^*012$ is a g-1-flip;
    \item [2:] $\flp{2(01)^+021}2=120(10)^+2$ is a g-1-flip;
    \item [2:] $\flp{21}2(01)^+02=12(01)^+02$ is a g-1-flip.
    \end{enumerate}

\item [Type V, odd:] $(10)^+212$, has possible parents starting with:
    \begin{enumerate}
    \item [0:] $(01)^+0212$ is of type V;
    \item [2:] $\flp{2(01)^+21}2=12(10)^+2$ is a g-1-flip;
    \item [2:] $\flp{21}2(01)^i2=12(01)^i2$ is a g-1-flip unless $i=1$, but then the parent is $212012=Y_3$.
    \end{enumerate}


\item [Type VI, $1(20)^+1(20)^*2$:] has possible parents starting with:
    \begin{enumerate}
    \item [0:] $(02)^i10(20)^j1(20)^k2$ with $i>0$:
        \begin{description}
        \item [If] $i>1$: $\flp{02}02(02)^{i-2}10(20)^j1(20)^k2=2(02)^{i-1}10(20)^j1(20)^k2$ is a g-1-flip;
        \item [If] $i=1$, $j>0$: $\flp{02102}0(20)^{j-1}1(20)^k2=201(20)^j1(20)^k2$ is a g-1-flip;
        \item [If] $i=1$, $j=0$, $k>0$: $\flp{021012}0(20)^{k-1}2=2101(20)^k2$ is a g-1-flip;
        \item [If] $i=1$, $j=k=0$: $021012=Y_2$;
        \end{description}
    \item [0:] $\flp{(02)^+1}(02)^+1(02)^+=1(20)^+21(02)^+$ is a g-1-flip;
    \item [2:] $\flp{2(02)^*1(20)^+1(20)^*}2=(02)^*1(02)^+1(20)^*2$ is a g-1-flip unless this last string is $10212$ (type V)
    or $0210212$ (type VI), but then the parent is $212012=Y_3$ or $21201202=Y_{28}$ respectively;
    \item [2:] $\flp{2(02)^*1(02)^+1(20)^+}2=(02)^+1(20)^+1(20)^*2$ is a g-1-flip;
    \item [2:] $\flp{2(02)^*10}2(02)^*12=01(20)^*212$ is a g-1-flip unless there is no second 0, but then the
    parent is $210212=Y_1$.
    \end{enumerate}

\item [Type VI, $0(21)^+0(21)^*2$:] has possible parents starting with:
    \begin{enumerate}
    \item [1:] $(12)^i01(21)^j0(21)^k2$ with $i>0$:
        \begin{description}
        \item [If] $i>1$: $\flp{12}12(12)^{i-2}01(21)^j0(21)^k2=2(12)^{i-1}01(21)^j0(21)^k2$ is a g-1-flip;
        \item [If] $i=1$, $j>0$: $\flp{12012}1(21)^{j-1}0(21)^k2=210(21)^j0(21)^k2$ is a g-1-flip;
        \item [If] $i=1$, $j=0$, $k>0$: $\flp{120102}1(21)^{k-1}2=2010(21)^k2$ is a g-1-flip;
        \item [If] $i=1$, $j=k=0$: $120102=Y_4$;
        \end{description}
    \item [1:] $\flp{(12)^+0}(12)^+0(12)^+=0(21)^+20(12)^+$ is a g-1-flip;
    \item [2:] $\flp{2(12)^*0(21)^+0(21)^*}2=(12)^*0(12)^+0(21)^*2$ is a g-1-flip unless this last string is $1201202$ (type VI), but
    then the parent is $20210212=Y_{29}$;
    \item [2:] $\flp{2(12)^*0(12)^+0(21)^+}2=(12)^+0(21)^+0(21)^*2$ is a g-1-flip;
    \item [2:] $\flp{2(12)^*01}2(12)^*02=10(21)^*202$ is a g-1-flip unless this last string is $10202$ (type III), but then
    the parent is $201202=Y_5$.
    \end{enumerate}

\item [Type VII:] $1(02)^+1(02)^+$ has possible parents starting with:
    \begin{enumerate}
    \item [0:] $\flp{0(20)^*1(02)^+1}(02)^+=1(20)^+1(02)^+$ is a g-1-flip;
    \item [0:] $\flp{0(20)^*12}0(20)^*1(02)^+=210(20)^*1(02)^+$ is a g-1-flip;
    \item [2:] $\flp{(20)^+12(02)^*10}2(02)^*=01(20)^*21(02)^+$ is a g-1-flip;
    \item [2:] $\flp{(20)^+1}(20)^+12(02)^*=1(02)^+012(02)^*$ is a g-1-flip.
    \end{enumerate}

\item [Type VIII:] $1(02)^+12$ has possible parents starting with:
    \begin{enumerate}
    \item [0:] $0(20)^i1(02)^j12$ with $j>0$:
        \begin{description}
        \item [If] $i>0$: $\flp{02}0(20)^{i-1}1(02)^j12=(20)^i1(02)^j12$ is a g-1-flip;
        \item [If] $i=0$, $j>1$: $\flp{0102}(02)^{j-1}12=201(02)^{j-1}12$ is a g-1-flip;
        \item [If] $i=0$, $j=1$: $010212$ is of type V;
        \end{description}
    \item [2:] $\flp{(20)^+12(02)^*1}2=1(20)^*21(02)^+$ is a g-1-flip;
    \item [2:] $21(20)^i12$ with $i>0$:
        \begin{description}
        \item [If] $i=1$: $212012=Y_3$;
        \item [If] $i>1$: $\flp{2120}(20)^{i-1}12=021(20)^{i-1}12$ is a g-1-flip. \QeD
        \end{description}
    \end{enumerate}
\end{description}
\end{proof}

\begin{theorem}
\label{sorting:thm:alg} There exists a polynomial time algorithm for optimally sorting ternary strings.
\end{theorem}
\begin{proof}
Follows rather easily from Theorem \ref{thm:diamref}.
\QeD
\end{proof}
Finally, in light of the fact that the complexity of the sorting problem on quaternary (and higher) strings remains open, the following
serves as an intermediate result:

\begin{theorem}
\label{sorting:thm:ptas}
For every fixed $k$ there is a PTAS for sorting $k$-ary strings.
\end{theorem}

\begin{proof}
The proof is very similar to the proof of Theorem \ref{grouping:thm:ptas}. We assume that $k\geq 4$.
Let $N = (3k-2)/\epsilon + k$. Let $s$, the string we wish to sort, be of length $n$.
We distinguish two cases. (In both cases it is useful to note that $d_{\rm s}(s) \leq 2n$ because
we can always bring the greatest symbol not yet in its final position to the front and then to its
correct position.)\\
\\ 
\textbf{Case 1.} If $n \geq N$, we first
group the string using the ``greedy'' algorithm from
the proof of Lemma \ref{grouping:lem:greedy}, which yields
a permutation on $k$ symbols. This permutation can then be easily sorted with at most $2k$ flips.
Thus the total number of flips, denoted by $d^{G}_{\rm s}(s)$, is at most
$(n-2) + 2k$. This together with the grouping lower bound of Lemma~\ref{grouping:lem:lowerbound}
of $n-k$ on $d_{\rm s}(s)$ yields
$d^{G}_{\rm s}(s) \leq d_{\rm s}(s) + (3k-2) \leq (1+\epsilon)d_{\rm s}(s)$.\\
\\
\textbf{Case 2.} If $n < N$ we apply brute force by selecting the shortest shorting sequence from
among all length-$2n$ sequences of flips; there are at most $n^{2n}$ such sequences.
Given that $d_{\rm s}(s) \leq 2n$ this is guaranteed to give an optimal solution. The running time
in this case is bounded by a constant. \QeD
\end{proof}

\section{Prefix reversal diameter}
\label{sec:diam}

Let $S(n,k)$ be the set of fully $k$-ary strings of length $n$.
We define $\delta(n,k)$ as the largest value of $d(s, t)$ ranging
over all compatible
$s,t \in S(n,k)$.

\begin{theorem}
\label{thm:bidam}
For all $n \geq 2$, $\delta(n, 2) = n-1$.
\end{theorem}
\begin{proof} To prove $\delta(n,2) \geq n-1$, consider compatible $s,t\in S(n,2)$ with
$s=(10)^{n/2}$ in case $n$ even and $s=0(10)^{(n-1)/2}$ in case $n$ odd and in both cases
$t=I(s)$ i.e.
$t$ is the sorted version of $s$. By Theorem \ref{sorting:thm:bin}, $d(s,t) \geq n-1$.

The proof that $\delta(n,2) \leq n-1$, for all $n \geq 2$
is by induction on $n$. The lemma is trivially true for $n=2$.
Consider two compatible binary strings of length $n$:
$s=s_1s_2\ldots s_n$ and $t=t_1t_2\ldots t_n$.
If $s_n=t_n$ then by induction $d(s,t)\leq n-2$.
Thus, suppose
(wlog) $s_n=0$ and $t_n=1$. If $t_1=0$ then $f^{(n)}t$ and $s$ both end
with a 0, and using induction and symmetry $d(s,t)\leq 1 + d(f^{(n)}t,s) \leq n-2+1=n-1$.
An analogous argument holds if $s_1=1$.

Remains the case $s_1= s_n=0$ and $t_1 = t_n=1$.
First, suppose $t_{n-1}=0$. Since $s$ and $t$ are compatible,
there must exist index $i$ such that $s_i=0$ and
$s_{i+1}=1$.
Hence, $f^{(n)}(f^{(i+1)}(s))$ ends with $01$ like $t$ and by induction
$d(s,t)\leq 2+d(f^{(n)}(f^{(i+1)}(s)),t) = 2+n-3$.
Analogously, we resolve the case $s_{n-1}=1$.

Finally, suppose $s = 0...00$ and $t = 1...11$.
If $s$ contains 11 as a substring, then flipping that 11 (in the same manner as above)
to the back of $s$ using 2 flips, gives two
strings that both end in 11. Alternatively,
if $s$ does not contain 11 as a substring then $s$ has
at least two more 0's than 1's, which implies that $t$
must contain 00 as a substring. In that case two prefix reversals on
$t$ suffice to create two strings that both end
with $00$. In both cases, the induction hypothesis gives the required bound.
\QeD
\end{proof}

Note that, trivially, $d(s,t) \leq 2n$
for all compatible $s,t\in S(n,k)$, for all $k$, because two prefix reversals
always suffice to increase the maximal common suffix between $s$ and $t$ by at least 1.
The following tighter bound gives the best bound known on the diameter of ternary strings.
\begin{lemma}
\label{lem:terdamub}
For any two compatible  $s,t \in S(n,k)$, for any $k$, let $a$ be the most frequent symbol in $s$ and $\alpha$ its multiplicity. Then $d(s,t) \leq 2(n-\alpha)$.
\end{lemma}
\begin{proof} We prove the lemma, by induction on $n$.
The lemma is trivially true for $n=2$.
Consider $s,t\in S(n,k)$. If $s_n=t_n=a$ then $s_1s_2\ldots s_{n-1}$ and $t_1t_2\ldots t_{n-1}$ are
compatible
length-$(n-1)$ strings where the most frequent symbol occurs at least $\alpha -1$ times. Thus, by induction $d(s,t) \leq 2((n-1)-(\alpha -1)) = 2(n-\alpha)$.
In case $s_n=t_n\neq a$ induction even gives $d(s,t) \leq 2((n-1)-(\alpha )) = 2(n- \alpha)-2$.
Thus, suppose $s_n\neq t_n$ implying wlog that $t_n = b \neq a$. Suppose $s_i=b$; after two flips
$s'=f^{(n)}(f^{(i)}(s))$ has $b$ at the end; $s'_n=t_n$. Moreover the length $n-1$ suffixes
of $s'$ and $t$ still contain $\alpha$ $a$'s. Hence by induction
$d(s,t) \leq 2+d(s',t) \leq 2 + 2((n-1)-\alpha)= 2(n-\alpha )$.
\QeD
\end{proof}

\begin{lemma}
\label{lem:terdamlb}
For all $n>3$, $n-1\leq \delta(n,3) \leq (4/3)n$.
\end{lemma}
\begin{proof}
Since in any ternary case $\alpha \geq \lceil n/3 \rceil$,
Lemma \ref{lem:terdamub} implies $\delta(n,3) \leq (4/3)n$.
To prove $\delta(n,3) \geq n-1$ we distinguish between $n$ is odd and $n$ is even.
For odd $n = 2h + 1$, let $s$ be $2(01)^{h}$, and for even $n = 2h$
let $s = 01(21)^{h-1}$. In both cases we let $t=I(s)$.
We observe that, in the even and in the odd case, $s2$ is a bad I-string and a bad IV-string, respectively,
in the sense of Definition \ref{sorting:def:badstrings}.
Thus, by Theorem \ref{thm:diamref} we have that $d(s,t) = d(s2,t2 ) = (n+1)-2 = n-1.$ (Here $s2$, respectively
$t2$, refers to the concatenation of $s$, respectively $t$, with an extra 2 symbol.)
\QeD
\end{proof}

Brute force enumeration has shown that, for $4 \leq n \leq 13$, $\delta(n,3) = n-1$. (Note that $\delta(3,3) = 3$ because
$d(021,012) = 3$.) Proving or disproving the conjecture
that $\delta(n,3) = n-1$ for $n>3$ remains an intriguing open problem\footnote{Interestingly, initial experiments with brute force enumeration have also shown that,
for $4 \leq n \leq 10$, $\delta(n,4) = n$, and for $5 \leq n \leq 9$, $\delta(n,5)=n$.}.

\section{Prefix reversal distance}
\label{sec:complexity}

We show that computing flip distance is NP-hard on binary strings.
We also point out, using a result from \cite{rsw}, that
computing flip distance on \emph{arbitrary} strings is polynomial-time reducible (in an
approximation-preserving sense) to computing it on binary strings.

\begin{theorem}
\label{theorem:nphard}
The problem of computing the prefix reversal distance of binary strings is NP-hard.
\end{theorem}
\begin{proof}
We prove NP-completeness of the corresponding decision problem:

\smallskip

\noindent
\emph{Name: }{\sc binary-PD} (2PD shortly)\\
\emph{Input: }Two compatible strings $s,t\in S(n,2)$, and a bound
$B \in \mathbb{Z}^{+}$.\\
\emph{Question: }Is $d(s, t) \leq B$?

\smallskip

\noindent
2PD$\in$NP, since a certificate for a positive answer consists of at most $B$ flips\footnote{Recall that for all
compatible strings $s, t\in S(n,2)$, trivially $d(s, t)\leq 2n$.}.
To show completeness we use a reduction from
{\sc 3-Partition} \cite{gjbook} (cf. \cite{christieirving} and \cite{rsw}).

\smallskip

\noindent
\emph{Name: }{\sc 3-Partition} (3P shortly)\\
\emph{Input: } A set $A=\{ a_1, a_2, ..., a_{3k}\} $ and
a number $N \in \mathbb{Z}^{+}$. Element $a_i$ has size $r(a_i) \in \mathbb{Z}^{+}$ satisfying $N/4 < r(a_i) < N/2$,
$i=1,\ldots,3k$, and $\sum_{i=1}^{3k}r(a_i) = kN$.\\
\emph{Question}: Can $A$ be partitioned into $k$ disjoint triplet sets $A_1$, $A_2$, ...,
$A_k$ such that $\sum_{a\in A_j}r(a) = N$, $j=1,\ldots,k$? \\

\noindent
Given instance $I = (A, N, r)$ of 3P, we create an instance of 2PD by setting $B=6k$
and building two compatible binary strings $s$ and $t$:
\[
s = \bigg ( \prod_{1 \leq i \leq 3k} 0001^{r(a_i)} \bigg )000
\ \ \ \ \ \ \ t = 0^{3(3k+1) - k}(01^{N})^{k}
\]
This construction is clearly polynomial in a unary encoding of the 3P instance; we use
the {\em strong} NP-hardness of 3P \cite{gjbook}. We claim that
$I = (A,N,r)$ is a positive instance of 3P $\Leftrightarrow$ $d(s, t) \leq 6k$.

\smallskip

\noindent
$\Rightarrow$)
Let $a_{ij}$ denote the $j$th element from triples $A_i$ (in arbitrary order), $j=1,2,3$, $i=1,\ldots,k$, and
let us abuse its name also to denote the corresponding 1-block of length $r(a_{ij})$ in $s$.

That $s$ can be transformed to $t$ in $6k$ flips follows directly from the
correctness of the following claim for $h=k$.

\smallskip

\noindent
\textbf{Claim.} For $0 < h \leq k$, $s$ can be transformed into a string $\psi_h = \alpha_h \omega_h$ in $h$ phases,
each consisting of 6 flips, where $\psi_h$ has the following specific properties:\\
(1) The suffix (i.e. $\omega_h$) is equal to $(01^N)^h$ and contains all $3h$ 1-blocks corresponding to the elements in $\cup_{j=1}^h
A_j$;\\
(2) The prefix (i.e. $\alpha_h$) contains the remaining $3(k-h)$ 1-blocks, each of them
flanked by 0-blocks of length at least 3, except possibly a 0-block of length 2 at its right end. (Given
that $\psi_h = \alpha_h \omega_h$ it follows that, in $\psi_h$, \emph{all} these remaining 1-blocks are
flanked by 0-blocks of length at least 3.)\\
\\
\emph{Proof.} The proof is by induction. First we transform $s$ into $\psi_1$ in 6
flips: flips 1 and 2 bring $a_{11}$ to the back, flips 3 and 4 bring $a_{12}$ to the back (just in front of $a_{11}$) and
flips 5 and 6 bring $a_{13}$ to the back (just in front of $a_{12}$). No 0-blocks are cut in this process, and
only 1-blocks $a_{11}, a_{12}$ and $a_{13}$ are affected (i.e. concatenated into a single length-$N$ 1-block).

Now, suppose by induction that after
$6(h-1)$ flips we have created $\psi_{h-1}$. The next 6 flips (which form phase $h$) work exclusively
on $\alpha_{h-1}$. Flips 1 and 2 bring $a_{h1}$ to the front and then to the back of $\alpha_{h-1}$; flips 3 and 4 bring $a_{h2}$ to the front and then to the back just in front of $a_{h1}$; flips 5 and 6 bring $a_{h3}$ to the front and then to the back just in front of $a_{h2}$. These 6 flips (which do not cut any 0-blocks
within $\alpha_{h-1}$)\footnote{Observe that, in terms of its action on the \emph{overall} string,
flip 2 of phase $h$ does
cut a 0-block, cutting $\alpha_{h-1}$ from $\omega_{h-1}$, creating the singleton 0-block in between
two length $N$ 1-blocks.}  thus transform
$\alpha_{h-1}$ into a string with $01^N$ at the suffix, which appended to $\omega_{h-1}$
gives a suffix equal to $\omega_{h}$. The only question is whether
the resulting overall string satisfies condition (2). The only obstacle to this is
the possible length-2 0-block at the end of $\alpha_{h-1}$. However, this block is not flipped in flip 1 of phase $h$, it is brought to the front in flip 2, and concatenated to another 0-block in flip 3, leaving the prefix string without a length-2 0-block. This completes the proof of the claim.\\

\noindent
$\Leftarrow$)
Suppose that $I$ is a negative instance of 3P. We show that $d(s, t) > 6k.$
Notice that if $I$ is not a positive instance then in any sequence
of flips taking $s$ to $t$ some flip must split a 1-block
i.e. $\flp{...1}1...$. Below we add this to a list of tasks that
any sequence of flips taking $s$ to $t$ must complete:

\smallskip

\noindent
(0) split at least one 1-block; \\
(1) reduce the number of 1-blocks by $2k$;\\
(2) bring a 1 symbol to the end of the string (because $t$ ends with a 1, but $s$ does not);\\
(3) increase the number of singleton 0-blocks by $k-1$;\\
(4) reduce the number of \emph{big} (i.e. of length at least 3) 0-blocks by $3k$.

\smallskip

\noindent
To prove that at least $6k+1$ flips are needed to complete tasks (0)-(4), we show
that flips which make progress towards completing one of the tasks can not effectively be used to
make progress on another task. From this it follows that at least
$1+2k+1+(k-1)+3k = 6k+1$ flips will be needed.

It is immediately clear that task (2), requiring a flip of a whole string,
cannot be combined with any of the other tasks in one flip.
Notice that any task(0)-flip (which is of the form $\flp{1...1}1...$ or of the form
$\flp{0...1}1...$) does not decrease the number of 1-blocks, while 0-blocks remain
unaffected. So such flips do not contribute to tasks (1)-(4).
Nor can any task(1)-flip  (which is always of the form $\flp{1...0}1...$) contribute
to any of the other tasks from the list. It is also not too difficult to verify that
it is not possible to reduce the number of big blocks by 2 or more in one flip.
%
%
However, some types of task(3)-flip \emph{can} at the same time also contribute to task (4), and
some other types of task(3)-flip can increase the number of singleton 0-blocks by two, effectively
contributing `twice' to task (3). Such flips we call (34)- and (33)- flips, respectively.
We will show that all (34)- and (33)-flips necessarily have to be succeeded by at least one flip that does not, in an
overall sense, help us with the completion of the tasks.

\noindent
Any (33)-flip is of the type \\
(33.1) $\flp{1...0}0...$ (where the 0s form a complete block)

\noindent
Any (34)-flip is of the type: \\
(34.1) $\flp{1...0}00...$ (where the 0s form a complete block)\\
(34.2) $\flp{1...00}0...$ (where the 0s form a complete block)\\
(34.3) $\flp{000...10}00...$ \\
We emphasize here that 00 is not considered to be a big 0-block.

After a flip of type (33.1), (34.1) or (34.3) we have a single 0 at the front. In such
a situation a task(1)- or task(2)- flip is not possible. We cannot perform a task(3)-flip
because flips of the form $\flp{01...}0...$ will destroy the initial singleton 0, and flips
of the form $\flp{01...}1...$ cannot create new singleton 0s. The only task(4)-flip
possible is $\flp{01...00}0...$ (where the second group of 0s forms a complete block) but this also
reduces the number of singleton 0-blocks by 1, meaning that an extra task(3)-flip
would then be needed. Termination is not an option (because $t$ does not begin
with $01$). A task(0)-flip of the form $\flp{01...1}1...$ is potentially possible but, as noted,
this increases the number of required task(1)-flips.\\
\\
After a flip of type (34.2) we are left with $001$ at the front. Again, a task(1)- or task(2)-
flip is not possible in this situation, and neither is termination. A task(3)-flip
is potentially possible but this brings a single 0 to the front, which (by the earlier
argument) cannot be followed by any useful flip. A task(4)-flip is not possible because, when
the string begins
with $001$, a task(4)-flip must necessarily split a 00-adjacency in some big 0-block, but this
simply creates a different big 0-block.
\QeD
\end{proof}

\noindent
For studying problems on
arbitrary strings, let $X$ and $Y$ be two compatible, length-$n$ strings, where we assume (wlog) that
each of the symbols from $X$ and $Y$ are drawn from the set $\{0,1,..., n-1\}$. We define $D(X,Y)$ as
the smallest number of flips required to transform $X$ to $Y$. The arity of the
strings $X$ and $Y$ does not need to be fixed, and symbols may be repeated. Hence,
sorting of a permutation by flips (MIN-SBPR), and the flip distance
problem over fixed arity strings, are both special cases of computing $D$. Given that
computing $D$ is a generalisation of computing distance $d$ of binary strings,
immediately implies that
it is NP-hard.
However, an approximation-preserving reduction in the \emph{other} direction is possible, meaning that
inapproximability results for one of the problems will be automatically inherited by the other.
\begin{theorem}
\label{thm:reduc}
Given two compatible strings $X$ and $Y$ of length $n$ with each symbol from $X$ and $Y$ drawn
from $\{0,1,..., n-1\}$, it is possible to compute in time polynomial in $n$
two binary strings $x$ and $y$ of length polynomial in $n$ such that $D(X,Y) = d(x,y)$.
\end{theorem}
As demonstrated shortly the above result follows directly from work by {\RSW}. A little background is necessary
to understand the context. In Theorem 8 of \cite{rsw} it is shown 
that sorting permutations by reversals is directly reducible to the reversal distance problem
on binary strings. It is later argued (in Theorem 11 of \cite{rsw}) that the same reduction
technique can be used to reduce the transposition distance problem on a 4-ary alphabet
to the transposition distance problem on a binary alphabet. The proof of Theorem 11 lacks detail but
personal communication with the authors \cite{rswP} has since clarified that the result is correct.
Furthermore, the reduction technique underpinning Theorems 8 and 11 from \cite{rsw} can be directly
applied to prove the present theorem. We show this by reproducing
the reduction technique (complete with clarification) in the context of prefix reversals.
We also use this opportunity to clarify the correctness of Theorem 11 from \cite{rsw}. The
following should thus be considered attributed to {\RSW}. \\
\\
\emph{Proof. } The strings $x$ and $y$ are constructed as follows:
\begin{align*}
x &= (10^{X_1 + 1}1)^{2n+1}...(10^{X_n + 1}1)^{2n+1}\\
y &= (10^{Y_1 + 1}1)^{2n+1}...(10^{Y_n + 1}1)^{2n+1}
\end{align*}
In the above encoding, each symbol $X_i$ is thus encoded as the fragment
$(10^{X_i + 1}1)^{2n+1}$, each fragment consisting of $2n+1$ subfragments. (This
also holds for each symbol in $Y$.) Note that a fragment is reversal invariant. To see that $d(x,y) \leq D(X,Y)$, observe
that - by mapping to prefix reversals that cut at the boundaries between fragments - any
sequence of $m$ prefix reversals taking $X$ to $Y$ can be trivially mapped to $m$ prefix
reversals which take $x$ to $y$.

The proof that $D(X,Y) \leq d(x,y)$ is more involved. Combining $d(x,y) \leq D(X,Y)$ with the
trivial fact that $D(X,Y) \leq 2n$ yields $d(x,y) \leq 2n$. Now, consider any shortest sequence of prefix reversals
taking $x$ to $y$. This sequence of prefix reversals will cut the string $x$ in at most $2n$ places.
A subfragment within $x$ is said to 
\emph{survive} iff it is not cut by any of these
prefix reversals. Now, construct a bipartite graph with vertex set $\{e_1, e_2, ..., e_n\} \cup \{f_1, f_2, ..., f_n\}$ and
add an edge $(e_i, f_j)$ iff some subfragment of the fragment corresponding to $X_{i}$ survives and ends
up in the fragment corresponding to $Y_{j}$. Observe that
within any set of $m$ fragments from $x$, strictly more than $(m-1)(2n+1)$ subfragments will
survive, and
hence at least $m$ fragments from $y$ will be required to absorb these surviving subfragments. Thus, by
Hall's Theorem, the graph has a perfect matching. For each edge $(e_i, f_j)$ of the perfect matching, pick a subfragment
from the fragment corresponding to $X_i$ that survives and ends up in the fragment corresponding to $Y_j$.
Considering the action of the flips only on these $n$ subfragments, we see that there exists a sequence
of $d(x,y)$ prefix reversals transforming the sequence of symbols in $X$ into the sequence of symbols in $Y$,
and thus $D(X,Y) \leq d(x,y)$. \QeD\\
\\
The correctness of Theorem 11 from \cite{rsw} follows by using the same reduction but encoding each fragment as $3n$ subfragments rather
than $2n+1$ subfragments. (The transposition distance between two compatible length-$n$ strings is strictly less than $n$, and a transposition cuts a 
string in at most 3 places.)  Indeed, it is easy to see that the reduction works for a whole family of string rearrangement operators, by
ensuring that the number of subfragments per fragment is sufficiently large. For example, consider a rearrangement operator
$op$, and let $u$ be some upper-bound on the number of places an $op$-operation can cut a string. Let $v$ be any upper
bound on the maximum value of $d_{op}(X,Y)$ ranging over all compatible length-$n$ strings $X,Y$. Encoding each fragment
with $uv+1$ subfragments is sufficient to generalise the above reduction.

%

\section{Open problems}
\label{sec:openproblems}

In this study we have unearthed many rich (and surprisingly difficult)
combinatorial questions which deserve further analysis. We discuss some of them here.
The main unifying, ``umbrella'' suggestion is that, to go beyond ad-hoc (and case-based)
proof techniques, it will be necessary to develop deeper, more structural insights
into the action of flips on strings over fixed size alphabets. \\
\\
\textbf{Grouping and sorting on higher arity alphabets. } We have shown how to
group and sort optimally binary and ternary strings, but characterisations and algorithms
for quaternary (and higher) alphabets have so far evaded us. As observed in Section
\ref{subsec:largeralphabets}, it \emph{seems} that for $k=4, 5$ and for sufficiently
long strings, the strings with grouping distance $n-2$ settle into some kind of pattern,
but this has not yet offered enough insights to allow either the development of
a characterisation or of an algorithm.
Related problems include: for all fixed $k$, are there polynomial algorithms to
optimally sort (optimally group) $k$-ary strings? Is grouping strictly easier than
sorting, in a complexity sense? How does grouping function under other operators
e.g. reversals, transpositions? An upper bound on the grouping transposition distance has been
presented in \cite{bridgehand}. \\
\\
\textbf{Diameter questions. } Proving or disproving that $\delta(n,3) = n-1$ for
$n>3$ remains the obvious open diameter question. Beyond that, diameter results
for quaternary and higher arity alphabets are needed. How does the diameter
$\delta(n,k)$ grow for increasing $k$? (At this point we conjecture that,
for sufficiently long strings, the diameter of 3-ary, 4-ary and 5-ary strings
is $n-1$, $n$, and $n$ respectively.)

 The suspicion also exists that, for
all $k$ and for all sufficiently long $n$, there exists a length-$n$ fully $k$-ary
string $s$ such that $d(s, I(s)) = \delta(n,k)$. In other words, the set of all
pairs of strings that are $\delta(n,k)$ flips apart includes some instances of
the sorting problem.
It should be noted however that, following empirical testing,
it is apparent that there are also very many pairs of strings $s,t$
with $s \neq I(t)$ and $t \neq I(s)$
that are $\delta(n,k)$ flips apart.

It also seems important to develop diameter results for subclasses of strings,
perhaps (as in \cite{rsw}) characterised by the frequency of their most frequent
symbol. It may be that such refined diameter results for $k$-ary alphabets
provide information that is important in determining $\delta(n,k+1)$.

Note finally that the diameter of
strings over fixed size alphabets, i.e. $\delta(n,k)$, is always bounded from above
by the diameter of permutations, $wc(n)$. This is because the distance problem
on two length-$n$, fixed size alphabet strings $s, t$ can easily be re-written as a
sorting problem on a length-$n$ permutation $\pi$, such that a sequence of prefix reversals
sorting the permutation also suffices to transform $s$ into $t$.
Indeed, because of this relabelling property, the flip distance between
two fixed size alphabet strings can be viewed as being equal to the minimum permutation sorting distance, ranging over all such
relabellings into a permutation $\pi$. Can this relationship
between the fixed size alphabet and permutation world be further specified and exploited?\\
\\
\textbf{Signed strings. } The problem of sorting signed permutations by flips
(the \emph{burnt} pancake flipping problem) is well known \cite{burnt} \cite{gatespap} \cite{heydarisud}, but in
this paper we have not yet attempted to analyse the action of flips
on signed, fixed size alphabet strings. Obviously, analogues of all the problems described
in this paper exist for signed strings.\\
\\
\textbf{Complexity/approximation. } In the presence of hardness results (e.g. Theorem
\ref{theorem:nphard}) it is interesting to explore the complexity of restricted
instances, and to develop algorithms with guaranteed approximation bounds. For example,
\cite{rsw} gives a PTAS for dense instances. The development
of approximation algorithms is also a useful intermediate strategy where the complexity
of a problem remains elusive. In particular, this requires the development of improved
lower bounds.

\section{Acknowledgements}

We thank Alex Scott for useful remarks made during the writing of this paper.


\newpage

\appendix


\setlength\oddsidemargin   {-15pt} \setlength\evensidemargin  {-15pt}

\begin{table}[h]
\begin{centering}
\begin{tabular}{|l|l|l|l|l|l|}
\hline $Y_{    1   }=  210212  $&$ P2:$ $\flp{12102}12  $&$ P3=Y_{13    }$&$    C4$ is of type I    &$  P5=Y_{15    }$&$
C6$ is of type VI  \\  \hline $Y_{    2   }=  021012  $&$ P2:$ $\flp{202101}2  $&$ P3=Y_{14    }$&$    C4$ is of type VI &$
P5:$ $\flp{10}12012  $&$ P6:$ $\flp{2101}202  $\\ \hline $Y_{    3   }=  212012  $&$ P2=Y_{11    }$&$    C3$ is of type VI &$
P4=Y_{8 }$&$    P5:$ $\flp{102}1212  $&$ C6$ is of type V    \\  \hline $Y_{    4   }=  120102  $&$ P2=Y_{9 }$&$ P3:$
$\flp{02101}02  $&$ C4$ is of type VI  &$  P5=Y_{10    }$&$    P6=Y_{12    }$\\    \hline $Y_{    5   }=  201202  $&$ P2:$
$\flp{02012}02  $&$ P3=Y_{7 }$&$    C4$ is of type III  &$  P5=Y_{6 }$&$    C6$ is of type VI  \\  \hline $Y_{    6 }=
0210202 $&$ P2:$ $\flp{20210}202 $&$ P3=Y_{35    }$&$    C4$ is of type II   &$  P5:$ $\flp{2012020}2 $&$ C6=Y_{5 }$\\&$ P7:$
$\flp{20}201202 $&$     $&$     $&$     $&$     $\\ \hline $Y_{    7   }=  1021202 $&$ P2=Y_{30    }$&$    P3:$ $\flp{20121}202
$&$ C4=Y_{5 }$&$    P5=Y_{28    }$&$    P6=Y_{23    }$\\&$  P7=Y_{32    }$&$        $&$     $&$     $&$ $\\ \hline
$Y_{    8   }=  0212012 $&$ P2=Y_{32    }$&$    P3=Y_{37    }$&$    P4=Y_{26    }$&$    C5=Y_{3 }$&$ P6=Y_{21 }$\\&$
P7:$ $\flp{210}21202 $&$     $&$     $&$     $&$     $\\ \hline $Y_{    9   }=  2120102 $&$ P2=Y_{33 }$&$ C3=Y_{4 }$&$
P4=Y_{23    }$&$    P5:$ $\flp{102}12102 $&$ P6=Y_{30    }$\\&$  C7$ is of type V    &$      $&$ $&$ $&$     $\\ \hline
$Y_{    10  }=  0102102 $&$ P2=Y_{24    }$&$    C3$ is of type VII   &$  P4:$ $\flp{2010210}2 $&$ P5:$ $\flp{120}10102 $&$
C6=Y_{4 }$\\&$  P7=Y_{36    }$&$        $&$     $&$     $&$     $\\ \hline $Y_{    11  }=  1212012 $&$ P2:$ $\flp{2121201}2
$&$ C3=Y_{3 }$&$    P4:$ $\flp{2121201}2 $&$ P5=Y_{22    }$&$    C6$ is of type I    \\&$    P7=Y_{18 }$&$ $&$     $&$ $&$
$\\ \hline $Y_{    12  }=  2010212 $&$ P2=Y_{25    }$&$    P3=Y_{17    }$&$ P4=Y_{31    }$&$ C5$ is of type V &$
P6=Y_{34    }$\\&$  C7=Y_{4 }$&$        $&$     $&$     $&$     $\\ \hline $Y_{    13  }= 0120212 $&$ P2:$ $\flp{101202}12
$&$ P3:$ $\flp{21020}212 $&$ C4=Y_{1 }$&$    P5=Y_{29    }$&$    P6:$ $\flp{1202}1012 $\\&$   P7:$ $\flp{21}202102 $&$ $&$     $&$
$&$     $\\ \hline $Y_{    14  }=  1201012 $&$ P2=Y_{27 }$&$    P3:$ $\flp{02101}012 $&$ C4=Y_{2 }$&$ P5:$ $\flp{01}021012 $&$
C6$ is of type V    \\&$    P7:$ $\flp{21010}212 $&$ $&$     $&$     $&$     $\\ \hline $Y_{    15  }= 1201212 $&$ P2:$
$\flp{21201}212 $&$ P3=Y_{20    }$&$    C4$ is of type I    &$  P5=Y_{18    }$&$ C6=Y_{1 }$\\&$  P7:$ $\flp{21}210212 $&$ $&$
$&$     $&$     $\\ \hline $Y_{    16  }= 2012012 $&$ P2:$ $\flp{02012}012 $&$ P3=Y_{21    }$&$
C4 $ is of type IV  &$  P5=Y_{19    }$&$    P6=Y_{17    }$\\&$ C7$ is of type VII   &$      $&$ $&$     $&$     $\\
\hline $Y_{    17  }=  10210212    $&$ P2=Y_{48    }$&$ P3:$ $\flp{201210}212    $&$ C4=Y_{12    }$&$
P5:$ $\flp{012}010212    $&$ P6:$ $\flp{201}201212    $\\&$   C7=Y_{16    }$&$ P8=Y_{47    }$&$        $&$     $&$     $\\
\hline $Y_{    18  }=  21021212    $&$ P2:$ $\flp{12102}1212    $&$ P3=Y_{56 }$&$    C4$ is of type I    &$  P5:$
$\flp{1201212}12    $&$ C6=Y_{15    }$\\&$  P7:$ $\flp{12}1201212    $&$ C8=Y_{11    }$&$ $&$     $&$     $\\ \hline $Y_{    19
}=  02102012    $&$ P2:$ $\flp{20210}2012    $&$ P3=Y_{57    }$&$    C4$ is of type VI  &$  P5=Y_{58    }$&$    C6=Y_{16
}$\\&$  P7:$ $\flp{1020}12012    $&$ P8:$ $\flp{210}201202    $&$     $&$     $&$ $\\ \hline $Y_{    20  }=  02101212    $&$
P2:$ $\flp{202101}212    $&$ P3:$ $\flp{12010}1212    $&$ C4=Y_{15    }$&$    P5:$
$\flp{10}1201212    $&$ P6:$ $\flp{2101}20212    $\\&$   P7:$ $\flp{12}1012012    $&$ P8:$ $\flp{21}2101202    $&$     $&$     $&$     $\\
\hline $Y_{    21  }=  10212012    $&$ P2=Y_{50    }$&$    P3:$ $\flp{20121}2012    $&$ C4=Y_{16    }$&$    P5=Y_{47 }$&$
P6=Y_{42    }$\\&$  C7=Y_{8 }$&$    P8:$ $\flp{210}212012    $&$     $&$     $&$     $\\ \hline $Y_{    22  }= 02121012 $&$
P2:$ $\flp{2021}21012    $&$ P3:$ $\flp{12012}1012    $&$ P4:$ $\flp{21}2021012    $&$ P5:$ $\flp{1212010}12    $&$ C6=Y_{11    }$\\&$
P7:$ $\flp{10}1212012    $&$ P8:$ $\flp{2101}21202    $&$     $&$     $&$     $\\ \hline $Y_{    23  }= 02120102    $&$ P2:$
$\flp{2021}20102    $&$ P3:$ $\flp{120120}102    $&$ P4:$ $\flp{21}2020102    $&$ C5=Y_{9 }$&$    P6=Y_{39 }$\\&$  C7=Y_{7 }$&$
P8=Y_{54    }$&$        $&$     $&$     $\\ \hline $Y_{    24  }=  10102102    $&$ P2:$ $\flp{0101021}02    $&$ C3=Y_{10
}$&$    P4:$ $\flp{0101021}02    $&$ P5:$ $\flp{20101210}2    $&$ C6 $ is of type III \\&$ P7:$ $\flp{012}010102    $&$ P8=Y_{60 }$&$
$&$     $&$     $\\ \hline $Y_{    25  }=  02010212    $&$ P2=Y_{51 }$&$    C3=Y_{12    }$&$    P4=Y_{40 }$&$    C5$
is of type VIII &$  P6:$ $\flp{20102021}2    $\\&$   P7=Y_{57    }$&$ P8=Y_{46    }$&$        $&$     $&$ $\\ \hline $Y_{
26  }=  21202012    $&$ P2=Y_{52    }$&$    C3$ is of type VI  &$  P4=Y_{41    }$&$    C5=Y_{8 }$&$ P6=Y_{45    }$\\&$
P7:$ $\flp{10202}1212    $&$ C8$ is of type VIII &$ $&$     $&$     $\\ \hline $Y_{    27  }=  21201012 $&$ P2:$
$\flp{1212010}12    $&$ C3=Y_{14    }$&$    P4=Y_{42 }$&$    P5:$ $\flp{102}121012    $&$ P6=Y_{50    }$\\&$  P7:$ $\flp{10}1021212
$&$ C8$ is of type V    &$      $&$     $&$ $\\ \hline $Y_{    28  }=  21201202    $&$ P2=Y_{53    }$&$ C3$ is of type
VI  &$  P4:$ $\flp{0212012}02    $&$ P5:$ $\flp{102}121202    $&$ C6=Y_{7 }$\\&$  P7=Y_{38    }$&$    C8$ is of type VI  &$ $&$
$&$     $\\ \hline $Y_{ 29  }=  20210212    $&$ P2=Y_{43    }$&$    C3$ is of type VI  &$ P4:$ $\flp{1202102}12 $&$ P5:$
$\flp{012}020212    $&$ C6=Y_{13    }$\\&$  P7=Y_{59    }$&$    C8$ is of type VI  &$ $&$     $&$     $\\ \hline $Y_{ 30 }=
01021202    $&$ P2:$ $\flp{10102}1202    $&$ C3=Y_{7 }$&$    P4=Y_{54    }$&$ P5:$ $\flp{120}101202    $&$ P6=Y_{46 }$\\&$
C7=Y_{9 }$&$    P8:$ $\flp{2021}20102    $&$     $&$     $&$     $\\ \hline $Y_{ 31  }=  01020212    $&$ P2=Y_{44 }$&$ C3$
is of type VIII &$  P4:$ $\flp{20102021}2    $&$ C5=Y_{12    }$&$    P6=Y_{51 }$\\&$  P7:$ $\flp{12020}1012 $&$ P8:$ $\flp{21}2020102
$&$     $&$     $&$     $\\ \hline $Y_{    32  }=  20212012    $&$ P2=Y_{45    }$&$    C3=Y_{8 }$&$
P4=Y_{55    }$&$    C5$ is of type VI  &$  P6=Y_{41    }$\\&$  P7:$ $\flp{102}120212 $&$ C8=Y_{7 }$&$ $&$     $&$ $\\
\hline
\end{tabular} \end{centering} \end{table} \begin{table}[h] \begin{centering} \begin{tabular}{|l|l|l|l|l|l|} \hline
$Y_{    33  }=  12120102    $&$ P2:$ $\flp{21212010}2    $&$ C3=Y_{9 }$&$    P4:$ $\flp{21212010}2    $&$ P5:$ $\flp{0212101}02 $&$
C6$ is of type VI  \\&$    P7:$ $\flp{01}0212102    $&$ P8:$ $\flp{2010}21212    $&$     $&$     $&$     $\\ \hline $Y_{ 34  }=
12010212    $&$ P2:$ $\flp{212010}212    $&$ P3:$ $\flp{02101}0212    $&$ C4$ is of type VI  &$  P5=Y_{48    }$&$ P6:$
$\flp{2010}21212    $\\&$   C7=Y_{12    }$&$    P8:$ $\flp{212010}212    $&$     $&$     $&$     $\\ \hline $Y_{    35  }=
12010202    $&$ P2=Y_{46    }$&$    P3:$ $\flp{02101}0202    $&$ C4=Y_{6 }$&$    P5=Y_{49    }$&$    P6=Y_{54    }$\\&$
P7:$ $\flp{02}0102102    $&$ P8=Y_{51    }$&$        $&$     $&$     $\\ \hline $Y_{    36  }=  20120102    $&$ P2:$
$\flp{02012}0102    $&$ P3=Y_{39    }$&$    C4 $ is of type III &$  P5:$ $\flp{021}020102    $&$ P6:$ $\flp{102}102102    $\\&$
P7=Y_{49    }$&$    C8=Y_{10    }$&$        $&$     $&$     $\\ \hline $Y_{    37  }=  12012012    $&$ P2=Y_{47 }$&$
P3:$ $\flp{021012}012    $&$ C4=Y_{8 }$&$    P5:$ $\flp{210}212012    $&$ P6:$ $\flp{021}021012    $\\&$   C7$ is of type VI &$ P8:$
$\flp{210}210212    $&$     $&$     $&$     $\\ \hline $Y_{    38  }=  021021202   $&$ P2:$ $\flp{20210}21202   $&$ P3=Y_{71
}$&$    C4 $ is of type II  &$  P5:$ $\flp{201}2021202   $&$ P6:$ $\flp{120}1201202   $\\&$   P7:$ $\flp{21201}20202 $&$ C8=Y_{28
}$&$    P9:$ $\flp{20}21201202   $&$     $&$     $\\ \hline $Y_{    39  }=  102120102   $&$ P2:$ $\flp{010212}0102 $&$ P3:$
$\flp{20121}20102   $&$ C4=Y_{36    }$&$    P5:$ $\flp{21}20120102   $&$ P6:$ $\flp{0212}010102   $\\&$ C7=Y_{23    }$&$ P8:$
$\flp{010212}0102   $&$ P9=Y_{70    }$&$        $&$     $\\ \hline $Y_{    40  }=  102010212   $&$ P2:$ $\flp{0102}010212 $&$
P3:$ $\flp{201201021}2   $&$ P4:$ $\flp{02}01010212   $&$ C5=Y_{25    }$&$    P6:$ $\flp{0102}010212   $\\&$ P7:$ $\flp{2010}201212 $&$ C8
$ is of type IV  &$  P9=Y_{68    }$&$        $&$     $\\ \hline $Y_{    41  }=  021202012 $&$ P2:$ $\flp{202120201}2 $&$
P3=Y_{72    }$&$    P4=Y_{67    }$&$    C5=Y_{26    }$&$    P6:$ $\flp{202120201}2   $\\&$ C7=Y_{32    }$&$    P8:$
$\flp{10202}12012   $&$ P9:$ $\flp{210}2021202   $&$     $&$     $\\ \hline $Y_{    42  }=  021201012 $&$ P2:$ $\flp{2021}201012 $&$
P3:$ $\flp{120120}1012   $&$ P4:$ $\flp{21}20201012   $&$ C5=Y_{27    }$&$    P6=Y_{63    }$\\&$ C7=Y_{21    }$&$    P8:$
$\flp{10102}12012   $&$ P9:$ $\flp{21010}21202   $&$     $&$     $\\ \hline $Y_{    43  }=  020210212 $&$ P2:$ $\flp{2020210}212 $&$
C3=Y_{29    }$&$    P4:$ $\flp{2020210}212   $&$ P5:$ $\flp{12020102}12   $&$ C6 $ is of type II
\\&$    P7:$ $\flp{20120}20212   $&$ P8=Y_{72    }$&$    P9:$ $\flp{21201}20202   $&$     $&$     $\\ \hline $Y_{    44  }=
101020212   $&$ P2:$ $\flp{010102}0212   $&$ C3=Y_{31    }$&$    P4:$ $\flp{010102}0212   $&$ P5:$ $\flp{2010120}212   $&$ P6:$
$\flp{02}01010212   $\\&$   P7:$ $\flp{20}20101212   $&$ C8 $ is of type IV  &$  P9:$ $\flp{21}20201012   $&$     $&$     $\\ \hline
$Y_{    45  }=  020212012   $&$ P2:$ $\flp{202021201}2   $&$ C3=Y_{32    }$&$    P4:$ $\flp{202021201}2   $&$ P5:$ $\flp{12020120}12
$&$ P6=Y_{67    }$\\&$  C7=Y_{26    }$&$    P8=Y_{62    }$&$    P9:$ $\flp{210}2120202   $&$     $&$     $\\ \hline $Y_{ 46
}=  212010202   $&$ P2:$ $\flp{12120}10202   $&$ C3=Y_{35    }$&$    P4:$ $\flp{021201}0202   $&$ P5:$ $\flp{102}1210202   $&$ P6:$
$\flp{01}02120202   $\\&$   C7=Y_{30    }$&$    P8:$ $\flp{0201}021202   $&$ C9=Y_{25    }$&$        $&$     $\\ \hline $Y_{ 47
}=  212012012   $&$ P2:$ $\flp{12120}12012   $&$ C3=Y_{37    }$&$    P4:$ $\flp{0212012}012   $&$ P5:$ $\flp{102}1212012 $&$ C6=Y_{21
}$\\&$  P7=Y_{61    }$&$    P8:$ $\flp{102}1021212   $&$ C9=Y_{17    }$&$        $&$     $\\ \hline $Y_{ 48  }= 010210212
$&$ P2:$ $\flp{10102}10212   $&$ C3=Y_{17    }$&$    P4:$ $\flp{2010210}212   $&$ P5:$ $\flp{120}1010212   $&$ C6=Y_{34 }$\\&$  P7:$
$\flp{201201021}2   $&$ P8:$ $\flp{120120}1012   $&$ P9:$ $\flp{21}20120102   $&$     $&$     $\\ \hline $Y_{    49  }= 010210202
$&$ P2=Y_{65    }$&$    C3$ is of type VII   &$  P4:$ $\flp{2010210}202   $&$ P5:$ $\flp{120}1010202 $&$ C6=Y_{35 }$\\&$  P7:$
$\flp{201201020}2   $&$ C8=Y_{36    }$&$    P9:$ $\flp{20}20120102   $&$     $&$     $\\ \hline $Y_{    50  }= 010212012   $&$
P2:$ $\flp{10102}12012   $&$ C3=Y_{21    }$&$    P4=Y_{70    }$&$    P5:$ $\flp{120}1012012 $&$ P6=Y_{68 }$\\&$  C7=Y_{27
}$&$    P8=Y_{63    }$&$    P9:$ $\flp{210}2120102   $&$     $&$     $\\ \hline $Y_{ 51  }=  202010212 $&$ P2=Y_{66    }$&$
C3=Y_{25    }$&$    P4=Y_{66    }$&$    P5=Y_{64    }$&$    P6:$ $\flp{0102}020212   $\\&$ C7=Y_{31    }$&$    P8:$
$\flp{120}1020212   $&$ C9=Y_{35    }$&$        $&$     $\\ \hline $Y_{ 52  }=  121202012   $&$ P2:$ $\flp{212120}2012   $&$
C3=Y_{26    }$&$    P4:$ $\flp{212120}2012   $&$ P5:$ $\flp{0212102}012   $&$
P6:$ $\flp{20}21212012   $\\&$   P7:$ $\flp{02}02121012   $&$ C8 $ is of type II  &$  P9:$ $\flp{21020}21212   $&$     $&$     $\\
\hline $Y_{    53  }=  121201202   $&$ P2:$ $\flp{2121201}202   $&$ C3=Y_{28    }$&$    P4:$ $\flp{2121201}202   $&$ P5:$
$\flp{02121012}02   $&$ C6 $ is of type II  \\&$    P7:$ $\flp{210}2121202   $&$ P8:$ $\flp{021}0212102   $&$ P9=Y_{69    }$&$ $&$
$\\ \hline $Y_{    54  }=  201021202   $&$ P2:$ $\flp{0201}021202   $&$ P3:$ $\flp{102102}1202   $&$ P4:$ $\flp{01}02021202 $&$
C5=Y_{30    }$&$    P6=Y_{71    }$\\&$  C7=Y_{35    }$&$    P8:$ $\flp{021201}0202   $&$ C9=Y_{23    }$&$        $&$ $\\
\hline $Y_{    55  }=  120212012   $&$ P2:$ $\flp{2120}212012   $&$ P3=Y_{61    }$&$    P4:$ $\flp{20}21212012   $&$ C5=Y_{32
}$&$    P6:$ $\flp{2120}212012   $\\&$   P7:$ $\flp{0212}021012   $&$ C8 $ is of type II  &$  P9:$ $\flp{210}2120212 $&$     $&$ $\\
\hline $Y_{    56  }=  012021212   $&$ P2:$ $\flp{101202}1212   $&$ P3:$ $\flp{21020}21212   $&$ C4=Y_{18 }$&$    P5=Y_{69 }$&$
P6:$ $\flp{1202}101212   $\\&$   P7:$ $\flp{21}20210212   $&$ P8:$ $\flp{12}12021012   $&$ P9:$ $\flp{21}21202102   $&$     $&$ $\\ \hline
$Y_{    57  }=  120102012   $&$ P2=Y_{68    }$&$    P3:$ $\flp{02101}02012   $&$ C4=Y_{19    }$&$    P5:$ $\flp{01}02102012
$&$ P6=Y_{70    }$\\&$  P7:$ $\flp{02}01021012   $&$ C8=Y_{25    }$&$    P9:$ $\flp{210}2010212   $&$     $&$ $\\ \hline $Y_{
58  }=  201202012   $&$ P2:$ $\flp{02012}02012   $&$ P3=Y_{62    }$&$ C4 $ is of type IV  &$  P5:$ $\flp{0210202}012   $&$
C6=Y_{19    }$\\&$  P7:$ $\flp{02}02102012   $&$ P8=Y_{64    }$&$    C9$ is of type VII   &$      $&$ $\\ \hline $Y_{    59
}=  120120212   $&$ P2:$ $\flp{21201}20212   $&$ P3:$ $\flp{021012}0212 $&$ C4 $ is of type II  &$  P5:$ $\flp{210}2120212   $&$ P6:$
$\flp{021}0210212   $\\&$   P7=Y_{69    }$&$    C8=Y_{29    }$&$ P9:$ $\flp{21}20210212   $&$     $&$ $\\ \hline $Y_{    60  }=
201201012   $&$ P2:$ $\flp{02012}01012   $&$ P3=Y_{63 }$&$    C4 $ is of type IV  &$  P5:$ $\flp{021}0201012   $&$ P6:$
$\flp{102}1021012   $\\&$   P7:$ $\flp{01}02102012   $&$ P8:$ $\flp{10102}10212   $&$ C9=Y_{24    }$&$ $&$     $\\ \hline $Y_{    61
}=  210212012   $&$ P2:$ $\flp{20210}212012  $&$ P3:$ $\flp{120102120}12  $&$ C4=Y_{55    }$&$ P5:$ $\flp{201}20212012  $&$ P6:$
$\flp{120}12012012  $\\&$   P7=Y_{76    }$&$ C8=Y_{47    }$&$    P9:$ $\flp{102}12012012  $&$ P10:$ $\flp{210}21201202 $&$     $\\
\hline
\end{tabular} \end{centering} \end{table} \begin{table}[h] \begin{centering} \begin{tabular}{|l|l|l|l|l|l|} \hline
$Y_{    62  }=  1021202012  $&$ P2:$ $\flp{01021202}012  $&$ P3:$ $\flp{20121}202012  $&$ C4=Y_{58    }$&$    P5=Y_{76    }$&$
P6=Y_{74    }$\\&$  P7:$ $\flp{20}212012012  $&$ P8:$ $\flp{02}021201012  $&$ C9=Y_{45    }$&$    P10:$ $\flp{210}20212012 $&$ $\\
\hline $Y_{    63  }=  1021201012  $&$ P2:$ $\flp{010212}01012  $&$ P3:$ $\flp{20121}201012  $&$ C4=Y_{60    }$&$    P5:$
$\flp{21}201201012  $&$ P6:$ $\flp{0212}0101012  $\\&$   C7=Y_{42    }$&$    P8:$ $\flp{010212}01012  $&$ C9=Y_{50    }$&$    P10:$
$\flp{21010}212012 $&$     $\\ \hline $Y_{    64  }=  1020210212  $&$ P2:$ $\flp{0102}0210212  $&$ P3:$ $\flp{20120}210212  $&$ P4:$
$\flp{02}010210212  $&$ P5:$ $\flp{20}201210212  $&$ C6=Y_{51    }$\\&$  P7:$ $\flp{012}02010212  $&$ P8:$ $\flp{201}20201212  $&$
C9=Y_{58    }$&$    P10=Y_{76   }$&$        $\\ \hline $Y_{    65  }=  1010210202  $&$ P2:$ $\flp{0101021}0202  $&$
C3=Y_{49    }$&$    P4:$ $\flp{0101021}0202  $&$ P5:$ $\flp{20101210}202  $&$ C6 $ is of type III \\&$    P7:$ $\flp{012}01010202 $&$
P8:$ $\flp{201}20101202  $&$ P9:$ $\flp{02}012010102  $&$ P10:$ $\flp{20}201201012 $&$     $\\ \hline $Y_{    66  }=  0202010212 $&$
P2:$ $\flp{20202010}212  $&$ C3=Y_{51    }$&$    P4:$ $\flp{20202010}212  $&$ C5=Y_{51    }$&$    P6=Y_{73    }$\\&$  C7$ is of
type VIII &$  P8:$ $\flp{2010}2020212  $&$ P9:$ $\flp{120}10202012  $&$ P10:$ $\flp{21}201020202 $&$     $\\ \hline $Y_{    67  }=
2120202012  $&$ P2=Y_{77    }$&$    C3$ is of type VI  &$  P4:$ $\flp{021202}02012  $&$ C5=Y_{41    }$&$    P6:$
$\flp{02}021202012  $\\&$   C7=Y_{45    }$&$    P8:$ $\flp{02}020212012  $&$ P9:$ $\flp{1020202}1212  $&$ C10$ is of type VIII    &$
$\\ \hline $Y_{    68  }=  2120102012  $&$ P2:$ $\flp{12120}102012  $&$ C3=Y_{57    }$&$    P4=Y_{74    }$&$    P5:$
$\flp{102}12102012  $&$ P6:$ $\flp{01021202}012  $\\&$   C7=Y_{50    }$&$    P8:$ $\flp{0201}0212012  $&$ P9:$ $\flp{1020}1021212  $&$
C10=Y_{40   }$&$        $\\ \hline $Y_{    69  }=  2021021212  $&$ P2=Y_{75    }$&$    C3$ is of type VI  &$  P4:$
$\flp{1202102}1212  $&$ P5:$ $\flp{012}02021212  $&$ C6=Y_{56    }$\\&$  P7:$ $\flp{120120212}12  $&$ C8=Y_{59    }$&$    P9:$
$\flp{12}120120212  $&$ C10=Y_{53   }$&$        $\\ \hline $Y_{    70  }=  2010212012  $&$ P2:$ $\flp{0201}0212012  $&$ P3:$
$\flp{102102}12012  $&$ P4:$ $\flp{01}020212012  $&$ C5=Y_{50    }$&$    P6:$ $\flp{120102120}12  $\\&$   C7=Y_{57    }$&$ P8=Y_{74
}$&$    P9:$ $\flp{102}12010212  $&$ C10=Y_{39   }$&$        $\\ \hline $Y_{    71  }=  1201021202  $&$ P2:$ $\flp{212010}21202
$&$ P3:$ $\flp{02101}021202  $&$ C4=Y_{38    }$&$    P5:$ $\flp{01}021021202  $&$ P6:$ $\flp{2010}2121202  $\\&$ C7=Y_{54    }$&$
P8:$ $\flp{212010}21202  $&$ P9:$ $\flp{0212}0102102  $&$ P10:$ $\flp{20}212010212 $&$     $\\ \hline $Y_{    72 }=  1201202012  $&$
P2=Y_{76    }$&$    P3:$ $\flp{021012}02012  $&$ C4=Y_{41    }$&$    P5:$ $\flp{210}21202012  $&$ P6:$ $\flp{021}02102012  $\\&$ P7:$
$\flp{20210}212012  $&$ P8:$ $\flp{02}021021012  $&$ C9=Y_{43    }$&$    P10:$ $\flp{210}20210212 $&$ $\\ \hline $Y_{    73  }=
10202010212 $&$ P2:$ $\flp{0102}02010212 $&$ P3:$ $\flp{20120}2010212 $&$ P4:$ $\flp{02}0102010212 $&$ P5:$ $\flp{20}2012010212 $&$ P6:$
$\flp{02}0201010212 $\\&$   C7=Y_{66    }$&$    P8:$ $\flp{0102}02010212 $&$ P9:$ $\flp{2010}20201212 $&$ C10 $ is of type IV &$
P11:$ $\flp{21}2010202012    $\\ \hline $Y_{    74  }=  02120102012 $&$ P2:$ $\flp{2021}20102012 $&$ P3:$ $\flp{120120}102012 $&$
P4:$ $\flp{21}2020102012 $&$ C5=Y_{68    }$&$    P6:$ $\flp{102}120102012 $\\&$   C7=Y_{62    }$&$ P8:$ $\flp{2010}21202012 $&$
C9=Y_{70    }$&$    P10:$ $\flp{1020}10212012    $&$ P11:$ $\flp{210}201021202    $\\ \hline $Y_{    75 }=  02021021212 $&$ P2:$
$\flp{2020210}21212 $&$ C3=Y_{69    }$&$    P4:$ $\flp{2020210}21212 $&$ P5:$ $\flp{12020102}1212 $&$ C6 $ is of type II  \\&$ P7:$
$\flp{201}202021212 $&$ P8:$ $\flp{120}120201212 $&$ P9:$ $\flp{21}2012020212 $&$ P10:$ $\flp{12120}1202012 $&$ P11:$ $\flp{21}2120120202 $\\
\hline $Y_{    76  }=  21201202012 $&$ P2:$ $\flp{12120}1202012 $&$ C3=Y_{72    }$&$ P4:$ $\flp{0212012}02012 $&$ P5:$
$\flp{102}121202012 $&$ C6=Y_{62    }$\\&$  P7:$ $\flp{021}021202012 $&$ C8=Y_{61    }$&$    P9:$ $\flp{02}0210212012 $&$ P10:$
$\flp{10202}1021212    $&$ C11=Y_{64   }$\\ \hline $Y_{    77  }=  12120202012 $&$ P2:$ $\flp{212120}202012 $&$ C3=Y_{67    }$&$
P4:$ $\flp{212120}202012 $&$ P5:$ $\flp{0212102}02012 $&$ P6:$ $\flp{20}2121202012 $\\&$ P7:$ $\flp{02}0212102012 $&$ P8:$ $\flp{20}2021212012
$&$ P9:$ $\flp{02}0202121012 $&$ C10 $ is of type II &$  P11:$ $\flp{210}202021212 $\\ \hline

\end{tabular}
\end{centering}
\caption{All strings of type IX (first column). For each string all parents and all
1-flips are listed. Each parent is either bad or a 1-flip to a good string is given.
For each string of type IX is also shown that each 1-flip leads to a bad string. Here $Pi$ denotes
the parent you get by doubling the $i$-th symbol and applying $p(i)$, $Ci$ denotes the string you get by applying the
1-flip $p(i-1)$. Note that if the $i$-th symbol is \emph{not} equal to the first symbol there is a parent
$Pi$ and if the $i$-th symbol \emph{is} equal to the first symbol there is a 1-flip possible, leading to
$Ci$.}
\label{sorting:veriftab}
\end{table}


\begin{thebibliography}{77} 

\bibitem{chen} X. Chen, J. Zheng, Z. Fu, P. Nan, Y. Zhong, S. Lonardi, T. Jiang,
Assignment of orthologous genes via genome rearrangement, IEEE/ACM Transactions
on Computational Biology and Bioinformatics, 2(4) (October-December 2005)

\bibitem{christieirving} D.A. Christie, R.W. Irving, Sorting strings by reversals and by transpositions,
SIAM J. on Discrete Math., 14(2) (2001), pp. 193-206

\bibitem{burnt} D.S. Cohen, M. Blum, On the problem of sorting burnt pancakes, Discrete Appl.
Math. 61(2) (1995), pp. 105-120

\bibitem{bridgehand} H. Eriksson, K. Eriksson, J. Karlander, L. Svensson, J. Wastlund,
Sorting a bridge hand, Discrete Math., 241 (2001), pp. 289-300

\bibitem{2approx} J. Fischer, S.W. Ginzinger, A 2-approximation algorithm for sorting by prefix reversals,
Proceedings of ESA 2005: 13th Annual European Symposium, LLNCS 3669 (2005), pp. 415-425

\bibitem{gj} M.R. Garey, D.S. Johnson, Complexity results for multiprocessor scheduling under
resource constraints, SIAM J. Comput., 4(4) (1975), pp. 397-411

\bibitem{gjbook} M.R. Garey, D.S. Johnson, Computers and Intractability: A Guide to the Theory
of NP-completeness, W. H. Freeman, San Francisco, CA, 1979

\bibitem{gatespap} W.H. Gates, C.H. Papadimitriou,
Bounds for sorting by prefix reversal, Discrete Math., 27 (1979), pp. 47-57

\bibitem{mcsp} A. Goldstein, P. Kolman, J. Zheng, Minimum Common String Partition problem: hardness
and approximations, Electron. J. Combin., 12, $\#$R50 (2005)

\bibitem{heydarisud} M.H. Heydari, I.H. Sudborough, On the diameter of the pancake network, J. Algorithms, 25 (1997) pp. 67-94

\bibitem{pnetwork} L. Morales, I.H. Sudborough, Comparing Star and Pancake Networks, The Essence
of Computation: Complexity, Analysis, Transformation, LNCS2566 (2002), pp. 18-36

\bibitem{rsw} A.J. Radcliffe, A.D. Scott and E.L. Wilmer, Reversals and transpositions over
finite alphabets, SIAM J. on Discrete Math., 19(1) (2005), pp. 224-244

\bibitem{rswP} A.D.Scott, Personal communication.

\end{thebibliography}
\end{document}